\documentclass[reqno,12pt]{amsart}
\usepackage{a4wide}
\usepackage{amsmath}       
\usepackage[foot]{amsaddr}
\usepackage{amsfonts}      
\usepackage{amssymb}       
\usepackage{amsthm}  
\usepackage{bm}            
\usepackage{enumitem}
\usepackage{graphicx}      
\usepackage{mathtools}     
\usepackage{xcolor}        
\usepackage[colorlinks,citecolor=green,linkcolor=red]{hyperref}
\usepackage[capitalize]{cleveref}      
\usepackage{dsfont}
\usepackage{tikz}
\usetikzlibrary{positioning,arrows.meta}
\usepackage{caption}
\usepackage{esint}

\usepackage{glossaries}
\expandafter\def\csname ver@etex.sty\endcsname{3000/12/31}

\usepackage{autonum}

\newcommand{\R}{\mathbb{R}}
\renewcommand{\H}{\mathcal{H}}
\renewcommand{\L}{\mathcal{L}}
\newcommand{\N}{\mathbb{N}}

\renewcommand{\d}{\,\mathrm{d}}
\newcommand{\rcd}{\mathrm{RCD}}

\newcommand{\Lip}{\mathrm{Lip}}
\newcommand{\Geo}{\mathrm{Geo}}
\newcommand{\lip}{\mathrm{lip}}
\newcommand{\supp}{\mathrm{supp}}
\newcommand{\diam}{\mathrm{diam}}

\renewcommand{\limsup}{\varlimsup}
\newcommand{\m}{\mathfrak{m}}

\usepackage[square,numbers]{natbib}

\numberwithin{equation}{section}
\numberwithin{equation}{section}

\newtheorem{theorem}{Theorem}[section]
\newtheorem{proposition}[theorem]{Proposition}
\AddToHook{env/proposition/begin}{\crefalias{theorem}{proposition}}

\AddToHook{env/conjecture/begin}{\crefalias{theorem}{conjecture}}

\AddToHook{env/question/begin}{\crefalias{theorem}{question}}

\newtheorem{lemma}[theorem]{Lemma}
\AddToHook{env/lemma/begin}{\crefalias{theorem}{lemma}}

\newtheorem{corollary}[theorem]{Corollary}
\AddToHook{env/corollary/begin}{\crefalias{theorem}{corollary}}

\theoremstyle{definition}
\newtheorem{definition}[theorem]{Definition}
\AddToHook{env/definition/begin}{\crefalias{theorem}{definition}}

\AddToHook{env/assumption/begin}{\crefalias{theorem}{assumption}}

\theoremstyle{remark}
\newtheorem{remark}[theorem]{Remark}
\AddToHook{env/remark/begin}{\crefalias{theorem}{remark}}

\newtheorem{example}[theorem]{Example}
\AddToHook{env/example/begin}{\crefalias{theorem}{example}}

\makeatletter
\renewcommand\subsubsection{\@startsection{subsubsection}{3}{\z@}%
	{.5\linespacing\@plus.7\linespacing}{-.5em}
	{\normalfont\bfseries}}

\let\emptyset\varnothing

\title[Lipschitz volume rigidity for non-collapsed $\rm RCD$ spaces with boundary]{Lipschitz volume rigidity for \\ non-collapsed $\rm RCD$ spaces with boundary}

\author{Zhenhao Li$^\dagger$}
\author{Enrico Pasqualetto$^\ddag$}
\author{Xin Qian$^\ddag$}
\address{$^\dagger$Faculty of Mathematics \\
	Bielefeld University \\
	Postfach 10 01 31 \\
	33501 Bielefeld \\
	Germany.}
\address{$^\ddag$Department of Mathematics and Statistics\\
	University of Jyv\"askyl\"a \\
	P.O.Box 35\\
	FI-40014\\
	Finland.}

\email{zhenhao.li@math.uni-bielefeld.de}
\email{enrico.e.pasqualetto@jyu.fi}
\email{xin.x.qian@jyu.fi}
\begin{document}
	
	\begin{abstract}
		We study Lipschitz-volume rigidity for non-collapsed ${\rm RCD}(K,N)$ spaces, possibly with boundary.
		We prove that a measure-preserving $1$-Lipschitz map between such spaces restricts to an isometry from the interior of the domain
		onto its image, where the latter is equipped with the intrinsic length distance. We then characterize when this interior rigidity
		extends to an isometric bijection between the whole spaces, in terms of the behavior of the map on the boundaries. Along the way,
		we establish a rigidity result of independent interest in a more general class of rectifiable metric spaces.
	\end{abstract}
	
	\keywords{Lipschitz-volume rigidity; non-collapsed ${\rm RCD}(K,N)$ spaces; rectifiable metric spaces; lip-to-Lip property.}
	\subjclass[2020]{53C23, 53C24}
	
	\maketitle
	
	\tableofcontents

	\section{Introduction}\label{sec:introduction}
	
	The Lipschitz-volume rigidity problem asks under what assumptions a $1$-Lipschitz map $f\colon X\to Y$ satisfying
	\begin{align}\label{eq:volume-equality}
		\H_X^n(X)=\H_Y^n(f(X))<\infty
	\end{align}
	is necessarily an isometry onto its image. The conclusion holds
	when both $X$ and $Y$ are closed connected Riemannian $n$-manifolds, but may fail for general length
	spaces, even if $f$ is a biLipschitz homeomorphism; see
	Example~\ref{example:shortcut-curve}.
	To include spaces of infinite volume, we consider the measure condition
	\[
	f_\#\H_X^n=\H_Y^n|_{f(X)}.
	\]
	For $1$-Lipschitz maps with $\H_X^n(X)<\infty$, this condition
	is equivalent to the volume equality \eqref{eq:volume-equality},
	and hence includes the finite-volume setting above; see
	Lemma~\ref{lemma:volume-preserving}.
	
	Li and Wang \cite{Li-Wang14} established Lipschitz-volume rigidity for
	non-collapsed limits of Riemannian manifolds with Ricci curvature
	bounded below. Under equality of the volumes of the domain and its
	image, they proved that a $1$-Lipschitz map is an isometry onto its image
	equipped with the intrinsic length distance. Li \cite{Li15} proved a
	corresponding result for Alexandrov spaces, allowing a finite disjoint
	union as the domain. In that setting, the map preserves the lengths of
	curves and is an intrinsic isometry on the interior.
	
	The problem has also been studied under different assumptions on the
	domain and target. Meier and Ntalampekos \cite{Meier-Ntalaepekos24}
	proved that every $1$-Lipschitz surjection from a closed metric surface
	onto a closed Riemannian surface of the same area is an isometry.
	They also obtained regularity results for maps between more general
	metric surfaces. Z\"ust \cite{Zust24} established rigidity for maps
	from integral current spaces to infinitesimally Euclidean Lipschitz
	manifolds, under suitable assumptions on the mass, boundary, and
	intrinsic metric.
	Related rigidity results for integral current spaces and metric manifolds
	were established in \cite{BassoCreutzSoultanis2023,DelNinPerales2023,Marti26}.
	
	A natural setting for this problem is provided by spaces satisfying
	the Riemannian curvature-dimension condition $\mathrm{RCD}(K,N)$.
	This condition combines a synthetic lower Ricci curvature bound with
	an infinitesimally Hilbertian differential structure; see
	\cite{LottVillani09,Sturm06-2,AGS14-RCD,Gigli15-RCD}.
	Following De Philippis and Gigli \cite{DePhilippis-Gigli18}, an
	$\mathrm{RCD}(K,N)$ space is called non-collapsed when its reference
	measure is the $N$-dimensional Hausdorff measure.
	Connell et al.\ \cite[Theorem~1.8]{CDNPSG24} extended the result
	of Li and Wang \cite{Li-Wang14} to non-collapsed $\mathrm{RCD}(K,N)$ spaces
	without boundary, for $N\geq3$.
	
	In this paper, we first prove a rigidity theorem for rectifiable
	metric spaces. We then establish intrinsic rigidity on the interior of
	non-collapsed $\mathrm{RCD}$ spaces with possibly non-empty boundary,
	and show that the map is an isometric bijection precisely when
	$f(\partial X)\subseteq\partial Y$.
	In particular, if the domain has empty boundary, global rigidity and
	surjectivity follow without any boundary assumption on the target.
	These rigidity results do not require compactness or surjectivity and
	allow infinite total volume. Finally, we obtain an
	almost-rigidity theorem for non-collapsed $\mathrm{RCD}$ spaces
	with uniform diameter and volume bounds.
	
	In the rectifiable setting, the measure identity and the $1$-Lipschitz bound
	force the tangential differential of $f$ to be an isometry almost everywhere.
	To recover distance preservation from this information, our first result
	uses the \emph{lip-to-Lip property}.
	For a Lipschitz function $\varphi\colon Z\to\R$, let
	$\Lip(\varphi)$ denote its global Lipschitz constant and set
	\[
	\lip(\varphi)(z)
	\coloneqq
	\limsup_{\substack{w\to z\\w\neq z}}
	\frac{|\varphi(w)-\varphi(z)|}{d_Z(w,z)},
	\]
	with value zero at isolated points.
	A metric measure space $(Z,d_Z,\mathfrak m)$ has the lip-to-Lip
	property if
	\[
	\Lip(\varphi)
	=
	\|\lip(\varphi)\|_{L^\infty(\mathfrak m)}
	\qquad
	\text{for every Lipschitz function }\varphi\colon Z\to\R.
	\]
	Thus an almost everywhere bound on the pointwise Lipschitz constant
	determines the global Lipschitz constant.
	As we explain in Appendix~\ref{app:lip-to-Lip}, this property holds
	for spaces satisfying the measure contraction property.
	It also follows from the Sobolev-to-Lipschitz property, thus in particular
	it holds for (possibly infinite-dimensional) $\mathrm{RCD}$ spaces.
	
	Recall that an $n$-rectifiable metric space is covered, up to an
	$\H^n$-null set, by countably many Borel sets admitting biLipschitz
	charts with values in $\R^n$.
	A set is porous if holes of a definite relative size occur at arbitrarily
	small scales near each of its points.
	
	\begin{theorem}[Lipschitz volume rigidity for rectifiable spaces]\label{thm:main_1}
		Let \((X,d_X)\) and \((Y,d_Y)\) be \(n\)-rectifiable, locally complete separable metric spaces such that
		\(\H^n_X\), \(\H^n_Y\) are boundedly finite, all porous sets in \(Y\) are
		\(\H^n\)-negligible, and \((Y,d_Y,\H^n_Y)\) has the lip-to-Lip property.
		Let \(f\colon X\to Y\) be an open, \(1\)-Lipschitz, locally biLipschitz map such that
		\(f_\#\H^n_X=\H^n_Y|_{f(X)}\). Then \(f\) is a local isometry and preserves the length of every
		rectifiable curve in $X$. If in addition \((X,d_X)\) and \((Y,d_Y)\) are complete,
		then \(f\colon X\to Y\) is an isometric bijection.
	\end{theorem}
	Here, a local isometry means that every point has a neighborhood on which the map preserves distances.
	In fact, we prove a slightly stronger result; see Theorem~\ref{thm:LVR-rectifiable}.
	The proof starts from the metric differentiation and area formula
	of Kirchheim \cite{Kirchheim94}.
	The assumption on porous sets allows us to relate tangential
	differentials to pointwise Lipschitz constants.
	We then apply the lip-to-Lip property to suitable functions
	constructed using local inverses of $f$.
	This gives the reverse distance inequality on small neighborhoods.
	
	The lip-to-Lip assumption cannot be omitted, even when both spaces
	are compact and the map is globally biLipschitz.
	The following example also explains why almost everywhere
	information alone does not determine the metric.
	
	\begin{example}[{\cite[Example~1.5]{Li15}}]
		\label{example:shortcut-curve}
		Let $(M,g)$ be a closed connected Riemannian manifold of dimension
		$n\geq2$, and let $\gamma$ be a nonconstant minimizing geodesic
		segment. Fix $\lambda\in(0,1)$.
		Define a length metric $d_\lambda$ by assigning weight $\lambda$
		to arclength along $\gamma$ and weight $1$ elsewhere.
		More precisely, for a Lipschitz curve $\sigma\colon[0,1]\to M$, set
		\[
		L_\lambda(\sigma)
		\coloneqq
		\int_0^1
		\bigl(1-(1-\lambda)\mathds{1}_\gamma(\sigma(t))\bigr)
		|\dot\sigma(t)|_g\,\mathrm dt,
		\]
		and let $d_\lambda$ be the infimum of these lengths over curves
		joining the given endpoints.
		Then
		\[
		\lambda d_g\leq d_\lambda\leq d_g.
		\]
		The two metrics agree locally away from $\gamma$, and $\gamma$
		has zero $n$-dimensional Hausdorff measure for both metrics.
		Consequently, the identity map
		\[
		f\colon(M,d_g)\longrightarrow(M,d_\lambda)
		\]
		is $1$-Lipschitz and globally biLipschitz, and satisfies
		$f_\#\H_{d_g}^n=\H_{d_\lambda}^n$.
		It is not a local isometry at interior points of $\gamma$.
		In particular, $(M,d_\lambda,\H_{d_\lambda}^n)$ does not have
		the lip-to-Lip property.
	\end{example}
	
	Further examples in Section~\ref{sec:LVR_rectifiable1} show the
	roles of local biLipschitz regularity, openness, and completeness.
	Whether the assumption on porous sets can be removed remains open.
	
	We now turn to non-collapsed $\mathrm{RCD}$ spaces.
	We write $\mathrm{ncRCD}(K,N)$ for an $\mathrm{RCD}(K,N)$ space
	whose reference measure is $\H^N$, and assume that this measure
	has full support.
	For such a space $Z$, let $\partial Z$ denote its boundary in the
	sense of \cite{DePhilippis-Gigli18}, and set
	\[
	Z^\circ\coloneqq Z\setminus\partial Z.
	\]
	For a subset $A\subset Z$, its intrinsic length distance is
	\[
	d_A^\ell(a,b)
	\coloneqq
	\inf\bigl\{
	L_Z(\sigma)\;\big|\;
	\sigma\colon[0,1]\to A\text{ is continuous},
	\ \sigma(0)=a,\ \sigma(1)=b
	\bigr\},
	\]
	where the infimum is $+\infty$ if there is no such curve.
	For the interior of an $\mathrm{ncRCD}$ space, $d_{Z^\circ}^\ell$
	coincides with the restriction of $d_Z$; see Lemma~\ref{lem:interior-ambient-metric}.
	
	\begin{theorem}[Lipschitz-volume rigidity for non-collapsed $\mathrm{RCD}$ spaces]
		\label{thm:main_2}
		Let $K\in\R$ and $N\in\N$, and let
		$(X,d_X,\H_X^N)$ and $(Y,d_Y,\H_Y^N)$ be
		$\mathrm{ncRCD}(K,N)$ spaces.
		Let $f\colon X\to Y$ be a $1$-Lipschitz map satisfying
		\(f_\#\H_X^N=\H_Y^N|_{f(X)}.\)
		Then $f$ preserves the length of every rectifiable curve in $X$, and $f|_{X^\circ}\colon X^\circ\to Y$ is an open map such that
		\[
		f|_{X^\circ}\colon
		(X^\circ,d_X)
		\longrightarrow
		\bigl(f(X^\circ),d_{f(X^\circ)}^\ell\bigr)
		\]
		is an isometry.
		Moreover, the following conditions are equivalent:
		\begin{enumerate}[label=\textup{(\roman*)}]
			\item $f(\partial X)\subseteq\partial Y$.
			\item $f(X^\circ)=Y^\circ$ and $f(\partial X)=\partial Y$.
			\item $f\colon(X,d_X)\to(Y,d_Y)$ is an isometric bijection.
		\end{enumerate}
	\end{theorem}
	
	The intrinsic distance in the conclusion is necessary.
	For example, the quotient map
	\[
	[0,\ell]\longrightarrow\R/\ell\mathbb Z
	\]
	is $1$-Lipschitz and satisfies the stated measure identity,
	but identifies the two boundary points.
	Its restriction to $(0,\ell)$ is an isometry for the intrinsic
	distance on the image, whereas distances measured in the circle
	can be smaller.
	Taking products with flat tori gives the same phenomenon in
	higher dimensions.
	The boundary condition in Theorem~\ref{thm:main_2} rules out
	this obstruction and also implies surjectivity.
	
	In particular, no assumption on the boundary of the target is
	needed when the domain has empty boundary.
	
	\begin{corollary}[Rigidity for spaces without boundary]
		\label{cor:main_boundaryless}
		Let $K\in\R$ and $N\in\N$, and let
		$(X,d_X,\H_X^N)$ and $(Y,d_Y,\H_Y^N)$ be
		$\mathrm{ncRCD}(K,N)$ spaces with $\partial X=\varnothing$.
		If a $1$-Lipschitz map $f\colon X\to Y$ satisfies
		\[
		f_\#\H_X^N=\H_Y^N|_{f(X)},
		\]
		then $f$ is an isometric bijection.
		In particular, $\partial Y=\varnothing$.
	\end{corollary}
	
	Theorem~\ref{thm:main_2} follows from Theorem~\ref{thm:length-isometry-RCD},
	Lemma~\ref{lemma:openness}, and Theorems~\ref{thm:LVR_RCD} and~\ref{thm:LVR_RCD_2}.
	Unlike Theorem~\ref{thm:main_1}, it does not assume that $f$
	is open or locally biLipschitz. We first use volume comparison and the
	rectifiable rigidity theorem to prove local isometry near regular points.
	Passing to tangent cones then shows that tangent maps preserve distances
	from the vertex. This implies that $f$ is asymptotically distance-preserving at every point of $X$,
	and hence preserves the lengths of all rectifiable curves, including those meeting the boundary.
	We next combine this property with a dimension-comparison argument,
	following the approach of Li and Wang \cite{Li-Wang14}, to prove openness on
	$X^\circ$. Volume preservation then gives injectivity there, and the same
	pointwise estimate controls curves under the inverse map and yields the
	intrinsic distance identity. Finally, the boundary condition allows us to
	extend this identity to an isometric bijection of the whole spaces.
	All arguments are local or take place on compact sets, so no global volume
	bound is needed.

	For Alexandrov spaces, the interior rigidity and length-isometry
	conclusions of Theorem~\ref{thm:main_2} follow from \cite{Li15}
	when $X$ has finite volume.
	Under the boundary condition $f(\partial X)\subseteq\partial Y$,
	our theorem yields an isometric bijection without the
	surjectivity assumption imposed in \cite[Corollary~0.1]{Li15}.
	
	The rigidity fails if the Lipschitz constant is allowed to exceed $1$,
	even arbitrarily slightly, as shown by the area-preserving maps
	$(x,y)\mapsto((1+\varepsilon)x,y/(1+\varepsilon))$ on $\R^2$,
	which are not isometries for any $\varepsilon>0$.
	
	Our final main result is an almost-rigidity theorem.
	It provides an $\mathrm{RCD}$ analogue of
	\cite[Corollary~0.1]{Li-Wang14}, while also allowing for spaces with non-empty boundary.
	
	\begin{theorem}[Almost rigidity]
		\label{thm:main_almost}
		Fix $K\in\R$, $N\in\N$, and $D,v>0$.
		For every $\varepsilon>0$, there exists $\delta=\delta(K,N,D,v,\varepsilon)>0$ such that the following holds.
		
		Let $(X,d_X,\H_X^N)$ and $(Y,d_Y,\H_Y^N)$ be ${\rm ncRCD}(K,N)$ spaces satisfying
		\begin{equation}
			\diam(X)\leq D,
			\qquad
			\H_Y^N(Y)\geq v.
		\end{equation}
		Suppose that a map $f\colon X\to Y$ satisfies
		\begin{equation}
			\begin{gathered}
				d_Y(f(x),f(x'))\leq d_X(x,x')+\delta \qquad\text{for all }x,x'\in X,\\
				\sup_{y\in Y}d_Y(y,f(X))\leq \delta,\qquad \H_X^N(X)\leq\H_Y^N(Y)+\delta,\\
				\H_X^{N-1}(\partial X)\leq\H_Y^{N-1}(\partial Y)+\delta.
			\end{gathered}
		\end{equation}
		Then $f$ is an $\varepsilon$-Gromov--Hausdorff approximation, i.e.
		\begin{equation}
			Y\subseteq B_\varepsilon(f(X)),\quad\left|d_Y(f(x),f(x'))-d_X(x,x')\right|\leq\varepsilon
			\qquad\text{for all }x,x'\in X.
		\end{equation}
	\end{theorem}
	
	No continuity assumption on $f$ is required in this statement.
	The proof, given in Theorem~\ref{thm:almost_rigidity}, combines a compactness argument with volume convergence \cite{DePhilippis-Gigli18} and boundary volume convergence \cite{Brue-Naber-Semola22}.
	A sequence of maps whose errors tend to zero yields a surjective
	$1$-Lipschitz map between limit spaces of equal volume.
	Boundary volume convergence also passes the boundary volume inequality to the limit.
	Proposition~\ref{prop:LVR_RCD_3} then identifies the limit map as an isometry.
	
	The paper is organized as follows.
	Section~2 develops the differential and measure-theoretic tools
	for rectifiable spaces and proves
	Theorem~\ref{thm:main_1}.
	Section~3 establishes the results for non-collapsed
	$\mathrm{RCD}$ spaces and the almost-rigidity theorem.
	Appendix~\ref{app:lip-to-Lip} studies the lip-to-Lip property,
	its relation to length spaces, and sufficient conditions for
	its validity.
	\subsubsection*{Acknowledgments}
	E.\ Pasqualetto and X.\ Qian were supported by the Research Council of Finland grant 362898.
	\section{Lipschitz volume rigidity for rectifiable metric spaces}
	
	Throughout this section, all ambient metric spaces are assumed to be locally complete and separable.\footnote{We adopt these standing assumptions for convenience. Separability is automatic in many settings considered here, for instance for rectifiable spaces.}
	These assumptions are inherited by open subsets.
	\subsection{Preliminaries on rectifiability}\label{sec:pre_rectifiable}
	For any \(n\in\N\), we denote by \(\L^n\) the Lebesgue measure in \(\R^n\). Given a Borel set \(A\subset\R^n\)
	and a Lipschitz map \(f\colon A\to\R^m\), we denote by
	\[
	{\rm ap}\d_x f\colon\R^n\to\R^m
	\]
	the \textbf{approximate differential} of \(f\) at \(x\), which
	exists for \(\L^n\)-a.e.\ \(x\in A\). More precisely, if \(\bar f\colon\R^n\to\R^m\) is any Lipschitz
	extension of \(f\) (whose existence follows by applying the McShane extension theorem coordinatewise),
	\(\text{ap}\d_x f\) coincides with the usual differential \(\d_x\bar f\) of \(\bar f\) at \(x\) for
	\(\L^n\)-a.e.\ \(x\in A\); the fact that \(\d_x\bar f\) exists for \(\L^n\)-a.e.\ \(x\in\R^n\) is
	guaranteed by the Rademacher theorem. See \cite{Federer69} for further details.
	\subsubsection*{Metric differential}
	Let \((X,d)\) be a metric space, \(A\subset\R^n\) a Borel set, and \(\phi\colon A\to X\) a Lipschitz map.
	It follows from the results of \cite{Kirchheim94} that for \(\L^n\)-a.e.\ \(x\in A\) there exists
	a unique seminorm \({\rm md}_x(\phi)\) on \(\R^n\), called the \textbf{metric differential} of \(\phi\)
	at \(x\), such that
	\[
	\lim_{A\ni y\to x}\frac{d(\phi(y),\phi(x))-{\rm md}_x(\phi)(y-x)}{|y-x|}=0.
	\]
	The \(\L^n\)-a.e.-defined field of seminorms \(A\ni x\mapsto{\rm md}_x(\phi)\) is measurable,
	meaning that \(A\ni x\mapsto{\rm md}_x(\phi)(v)\in\R\) is \(\L^n\)-measurable for every \(v\in\R^n\).
	Assuming in addition that \(\phi\colon A\to X\) is biLipschitz with its image, it holds that
	\({\rm md}_x(\phi)\) is a norm for \(\L^n\)-a.e.\ \(x\in A\).
	\subsubsection*{Rectifiable metric spaces}
	There is an extensive literature on rectifiable metric spaces. Although we rely primarily on \cite{Ambrosio-Kirchheim00},
	we adopt slightly different definitions in some places. For this reason, for the reader's convenience we provide a self-contained exposition.
	\medskip
	
	Let \((X,d)\) be a metric space. By an \textbf{\(n\)-chart} \((E,\varphi)\) on \(X\),
	for some \(n\in\N\), we mean a Borel set \(E\subset X\) with \(\H^n(E)>0\) together with a map \(\varphi\colon E\to\R^n\)
	that is biLipschitz with its image.
	In particular, if \(C\) denotes the biLipschitz constant of \(\varphi\colon E\to\varphi(E)\), then it holds that
	\[
	\frac{1}{C^n}\L^n|_{\varphi(E)}\leq\varphi_\#(\H^n|_E)\leq C^n\L^n|_{\varphi(E)}.
	\]
	We say that the metric space \((X,d)\) is
	\textbf{\(n\)-rectifiable} if there exists an \textbf{\(n\)-atlas} on \(X\), i.e.\ a countable collection
	\(\mathcal A=\{(E_i,\varphi_i):i\in\N\}\) of \(n\)-charts on \(X\) such that \(\{E_i:i\in\N\}\) is an
	\(\H^n\)-a.e.\ partition of \(X\), by which we mean that \(E_i\cap E_j=\varnothing\) for every \(i,j\in\N\)
	with \(i\neq j\), and \(\H^n(X\setminus\bigcup_{i\in\N}E_i)=0\). Given any \(i\in\N\), we define
	\[
	\|\cdot\|^{\mathcal A}_x\coloneqq{\rm md}_{\varphi_i(x)}(\varphi_i^{-1})\quad\text{ for }\H^n\text{-a.e.\ }x\in E_i.
	\]
	In this way, we obtain an \(\H^n\)-a.e.-defined measurable field of norms \(X\ni x\mapsto\|\cdot\|^{\mathcal A}_x\)
	on \(\R^n\). Whenever the chosen atlas \(\mathcal A\) is clear from the context, we simply write \(\|\cdot\|_x\)
	instead of \(\|\cdot\|^{\mathcal A}_x\). For every \(x\in X\) such that \(\|\cdot\|_x\) exists (thus, for
	\(\H^n\)-a.e.\ \(x\in X\)), we say that
	\[
	T_x X\coloneqq(\R^n,\|\cdot\|_x)
	\]
	is the \textbf{metric tangent} to \(X\) at \(x\).
	\begin{remark}{\rm
			Given two \(n\)-atlases \(\mathcal A=\{(E_i,\varphi_i):i\in\N\}\) and
			\(\tilde{\mathcal A}=\{(\tilde E_j,\tilde\varphi_j):j\in\N\}\) on \(X\), the tangent norms
			\(\|\cdot\|^{\mathcal A}_x\) and \(\|\cdot\|^{\tilde{\mathcal A}}_x\) associated with the two
			atlases are related by the following formula: for every \(i,j\in\N\) and \(\H^n\)-a.e.\ \(x\in E_i\cap\tilde E_j\),
			it holds that
			\[
			\|v\|^{\tilde{\mathcal A}}_x=\|\text{ap}\d_{\tilde\varphi_j(x)}(\varphi_i\circ\tilde\varphi_j^{-1})(v)\|^{\mathcal A}_x
			\quad\text{ for every }v\in\R^n.
			\]
			Note that \(\{E_i\cap\tilde E_j:i,j\in\N\}\) is an \(\H^n\)-a.e.\ partition of \(X\).
	}\end{remark}
	\subsubsection*{Tangential differentials}
	Let \((X,d_X)\) and \((Y,d_Y)\) be \(n\)-rectifiable metric spaces with atlases
	\(\mathcal A_X=\{(E_i,\varphi_i):i\in\N\}\) and \(\mathcal A_Y=\{(F_j,\psi_j):j\in\N\}\),
	respectively. Let \(f\colon X\to Y\) be a Lipschitz map. For any \(i,j\in\N\), we define
	\(G_f^{i,j}\coloneqq E_i\cap f^{-1}(F_j)\subset X\) and
	\[
	\tilde\d_x f\coloneqq\text{ap}\d_{\varphi_i(x)}f_{i,j}\colon T_x X\to\R^n
	\quad\text{ for }\H^n\text{-a.e.\ }x\in G_f^{i,j},
	\]
	where we set \(f_{i,j}\coloneqq\psi_j\circ f\circ\varphi_i^{-1}\colon\varphi_i(E_i\cap f^{-1}(F_j))\to\R^n\).
	Letting \(G_f\coloneqq\bigcup_{i,j\in\N}G_f^{i,j}\), we thus have an \(\H^n\)-a.e.\ map
	\(G_f\ni x\mapsto\tilde\d_x f\). However, for our purposes it would be desirable to have a differential
	defined as a map from \(T_x X\) to \(T_{f(x)}Y\). Since it can happen that \(T_{f(x)}Y\) does not exist
	for a set of points \(x\in G_f\) of positive \(\H^n\)-measure, we need to restrict
	the domain of definition of the differential to the possibly smaller set \(X_f\subset G_f\), which we
	define as \(X_f\coloneqq\bigcup_{i,j\in\N}X_f^{i,j}\), where we set
	\[
	X_f^{i,j}\coloneqq\big\{x\in G_f^{i,j}\;\big|\;\tilde\d_x f\colon T_x X\to\R^n\text{ is injective}\big\}
	\quad\text{ for every }i,j\in\N.
	\]
	The following result shows that the set \(X_f\) is well suited to our purposes.
	\begin{lemma}\label{lemma:reg_underf}
		Under the above assumptions, it holds that
		\[
		f_\#(\H^n|_{X_f^{i,j}})\ll\H^n\quad\text{ for every }i,j\in\N.
		\]
	\end{lemma}
	\begin{proof}
		Fix some Lipschitz extension \(\bar f_{i,j}\colon\R^n\to\R^n\) of \(f_{i,j}\). Given a Borel set \(N\subset\R^n\)
		such that \(\L^n(N)=0\), we denote \(S_N\coloneqq\varphi_i(X_f^{i,j})\cap f_{i,j}^{-1}(N)\). Since
		\(\tilde\d_{\varphi_i^{-1}(y)}f=\d_y\bar f_{i,j}\) for \(\L^n\)-a.e.\ point \(y\in S_N\), by applying the
		Euclidean area formula to \(\bar f_{i,j}\colon S_N\to N\) we obtain that
		\[
		\int_{S_N}|{\rm det}\tilde\d_{\varphi_i^{-1}(y)}f|\,\d\L^n(y)=\int_{S_N}|{\rm det}\d_y\bar f_{i,j}|\,\d\L^n(y)
		=\int_N\#(S_N\cap\bar f_{i,j}^{-1}(z))\,\d\L^n(z)=0.
		\]
		Given that \(|{\rm det}(\tilde\d_{\varphi_i^{-1}(y)}f)|>0\) for \(\L^n\)-a.e.\ \(y\in S_N\) by the injectivity
		of \(\tilde\d_{\varphi_i^{-1}(y)}f\), we deduce that \(\L^n(S_N)=0\), so that \(\H^n(\varphi_i^{-1}(S_N))=0\)
		and thus accordingly
		\[
		\H^n(X_f^{i,j}\cap(\psi_j\circ f)^{-1}(N))=\H^n(X_f^{i,j}\cap\varphi_i^{-1}(f_{i,j}^{-1}(N)))=\H^n(\varphi_i^{-1}(S_N))=0.
		\]
		This shows that \((\psi_j\circ f)_\#(\H^n|_{X_f^{i,j}})\ll\L^n\), whence it follows that \(f_\#(\H^n|_{X_f^{i,j}})\ll\H^n\).
	\end{proof}

	Taking Lemma \ref{lemma:reg_underf} into account, we get that \(\|\cdot\|_x\) and \(\|\cdot\|_{f(x)}\)
	exist for \(\H^n\)-a.e.\ \(x\in X_f^{i,j}\). Therefore, it makes sense to define the \textbf{tangential differential}
	\[
	\d_x f\colon T_x X\to T_{f(x)}Y\quad\text{ for }\H^n\text{-a.e.\ }x\in X_f
	\]
	by setting \(\d_x f\coloneqq\tilde\d_x f\) for every \(i,j\in\N\) and \(\H^n\)-a.e.\ \(x\in X_f^{i,j}\).
	Recall that \(\d_x f\) is linear.
	\medskip
	
	Fix a Lipschitz function \(\phi\colon Y\to\R\). In this case, the target space \(\R\) is tacitly assumed
	to be equipped with the atlas consisting solely of the identity chart \({\rm id}\colon\R\to\R\), and the metric
	tangents to \(\R\) can be identified with \(\R\) itself. Hence, we have that \(\H^n(Y\setminus G_\phi)=0\),
	and thus \(\tilde\d_y\phi\colon T_y Y\to\R\) exists for \(\H^n\)-a.e.\ \(y\in Y\). By applying Lemma
	\ref{lemma:reg_underf} and standard calculus rules, one can easily obtain the following chain rule.
	\begin{lemma}\label{lemma:chainrule}
		Under the above assumptions, both \(\tilde\d_x(\phi\circ f)\) and \(\tilde\d_{f(x)}\phi\) exist, and
		\[
		\tilde\d_x(\phi\circ f)=\tilde\d_{f(x)}\phi\circ\d_x f\quad\text{ for }\H^n\text{-a.e.\ }x\in X_f.
		\]
	\end{lemma}
	
	\subsubsection*{Jacobians and area formula}\label{sec:Jacobian}
	In accordance with \cite{Ambrosio-Kirchheim00}, the \textbf{\(n\)-Jacobian} \({\rm J}_n(s)\)
	of a seminorm \(s\) defined on \(\R^n\) is given by
	\[
	{\rm J}_n(s)\coloneqq\frac{\omega_n}{\L^n(\{x\in\R^n:s(x)\leq 1\})},
	\]
	whereas the \textbf{\(n\)-Jacobian} \({\rm J}_n(L)\) of a linear operator \(L\colon V\to W\)
	between two Banach spaces \(V\) and \(W\) with \({\rm dim}V=n\) is given by
	\[
	{\rm J}_n(L)\coloneqq\frac{\omega_n}{\H^n(\{v\in V:\|L(v)\|_W\leq 1\})};
	\]
	here, \(\omega_n\) denotes the Lebesgue measure of the unit ball in \(\R^n\). It follows from
	\cite[Lemma 6(ii)]{Kirchheim94} that for any norm \(\|\cdot\|\) on \(\R^n\) it holds that
	\begin{equation}\label{eq:Kirchh_ident}
		\H^n_{\|\cdot\|}={\rm J}_n(\|\cdot\|)\L^n.
	\end{equation}
	Note also that \({\rm J}_n(L)>0\) if and only if \(L\) is injective, thus accordingly
	\begin{equation}\label{eq:equiv_char_X_f}
		X_f=\{x\in G_f\;|\;{\rm J}_n(\tilde\d_x f)>0\}\quad\text{ up to }\H^n\text{-null sets.}
	\end{equation}
	\begin{lemma}\label{lemma:injectivePartition}
		There exists a sequence \((A_k)_{k\in\N}\) of pairwise disjoint compact subsets of \(X_f\) such that
		\(\H^n(X_f\setminus\bigcup_{k\in\N}A_k)=0\) and \(f|_{A_k}\) is injective for every \(k\in\N\).
	\end{lemma}
	\begin{proof}
		This is due to \cite[Lemma 2.9]{Bate-Takavc24} when $Y=\R^n$.
		For $Y$ being a general $n$-rectifiable space, we apply \cite[Lemma 2.9]{Bate-Takavc24} to each
		$f_j\colon X_f^{i,j}\to \psi_j(F_j)\subset \R^n$, where $f_j\coloneqq \psi_j\circ f$.
	\end{proof}

	\begin{definition}\label{def:J_n^x}
		We define
		\[
		{\rm J}_n^x f\coloneqq\left\{\begin{array}{ll}
			{\rm J}_n(\d_x f)\\
			0
		\end{array}\quad\begin{array}{ll}
			\text{ for }\H^n\text{-a.e.\ }x\in X_f,\\
			\text{ for }\H^n\text{-a.e.\ }x\in X\setminus X_f.
		\end{array}\right.
		\]
	\end{definition}

	It can be readily checked that the function \(X\ni x\mapsto{\rm J}_n^x f\in[0,+\infty)\) is \(\H^n\)-measurable.
	\begin{lemma}\label{lemma:Jacobians}
		For any \(i,j\in\N\), we have that
		\begin{equation}\label{eq:md&tandifferential}
			{\rm md}_{\varphi_i(x)}(f\circ\varphi_i^{-1})(v)=\|\d_ x f(v)\|_{f(x)}\quad\text{ for }\H^n\text{-a.e.\ }x\in X_f^{i,j}\text{ and every }v\in T_x X.
		\end{equation}
		In particular, we have that
		\begin{equation}\label{eqpf:Jacobians}
			{\rm J}_n^x f=\frac{{\rm J}_n({\rm md}_{\varphi_i(x)}(f\circ \varphi^{-1}_i))}{{\rm J}_n(\|\cdot\|_x)}
			\quad\text{ for }\H^n\text{-a.e.\ }x\in G_f^{i,j}.
		\end{equation}
	\end{lemma}
	\begin{proof}
		By unwrapping the various definitions, and using Lemma \ref{lemma:reg_underf} and the biLipschitzianity of
		the map \(\varphi_i\colon X_f^{i,j}\to\varphi_i(X_f^{i,j})\), for \(\H^n\)-a.e.\ point \(x\in X_f^{i,j}\) we obtain that
		\[\begin{split}
			d_Y(f(\tilde x),f(x))&={\rm md}_{\varphi_i(x)}(f\circ\varphi_i^{-1})(\varphi_i(\tilde x)-\varphi_i(x))+o(d_X(\tilde x,x)),\\
			d_Y(f(\tilde x),f(x))&=\|\psi_j(f(\tilde x))-\psi_j(f(x))\|_{f(x)}+o(d_Y(f(\tilde x),f(x))),\\
			\psi_j(f(\tilde x))-\psi_j(f(x))&=f_{i,j}(\varphi_i(\tilde x))-f_{i,j}(\varphi_i(x))=\d_x f(\varphi_i(\tilde x)-\varphi_i(x))+o(d_X(\tilde x,x)),
		\end{split}\]
		where \(\tilde x\in X_f^{i,j}\). Combining the above with the fact that \(x\) can be chosen so that \(\varphi_i(x)\) is a Lebesgue
		density point of \(\varphi_i(X_f^{i,j})\), we obtain \eqref{eq:md&tandifferential}. In particular, for \(\H^n\)-a.e.\ \(x\in X_f^{i,j}\)
		we have that
		\[\begin{split}
			{\rm J}_n^x f&=\frac{\omega_n}{\H^n(\{v\in T_x X:\|\d_x f(v)\|_{f(x)}\leq 1\})}
			=\frac{\omega_n}{{\rm J}_n(\|\cdot\|_x)\L^n(\{v\in T_x X:\|\d_x f(v)\|_{f(x)}\leq 1\})}\\
			&=\frac{\omega_n}{{\rm J}_n(\|\cdot\|_x)\L^n(\{v\in T_x X:{\rm md}_{\varphi_i(x)}(f\circ\varphi_i^{-1})(v)\leq 1\})}
			=\frac{{\rm J}_n({\rm md}_{\varphi_i(x)}(f\circ \varphi^{-1}_i))}{{\rm J}_n(\|\cdot\|_x)}
		\end{split}\]
		thanks to \eqref{eq:Kirchh_ident} and \eqref{eq:md&tandifferential}, thus proving \eqref{eqpf:Jacobians} in the case where
		\(x\in X_f^{i,j}\). To conclude, we consider the case where \(x\in G_f^{i,j}\setminus X_f^{i,j}\). For \(\H^n\)-a.e.\ such
		\(x\) we have that \({\rm ker}(\tilde\d_x f)\neq\{0\}\), so that we can find a sequence \(G_f^{i,j}\ni x_k\to x\) such that
		\[
		\tilde\d_x f(\varphi_i(x_k)-\varphi_i(x))=o(|\varphi_i(x_k)-\varphi_i(x)|)=o(d_X(x_k,x)).
		\]
		Arguing as above and using the fact that \(\psi_j\) is biLipschitz with its image, we obtain that
		\(d_Y(f(x_k),f(x))=o(d_X(x_k,x))\). Hence, the seminorm \({\rm md}_{\varphi_i(x)}(f\circ \varphi^{-1}_i)\)
		is not a norm, thus \({\rm J}_n({\rm md}_{\varphi_i(x)}(f\circ \varphi^{-1}_i))/{\rm J}_n(\|\cdot\|_x)=0={\rm J}_n^x f\).
		This completes the proof of \eqref{eqpf:Jacobians}.
	\end{proof}
	\begin{theorem}[Area formula]\label{thm:area_formula}
		Let \((X,d_X)\) and \((Y,d_Y)\) be \(n\)-rectifiable metric spaces with \(n\)-atlases \(\mathcal A_X\)
		and \(\mathcal A_Y\), respectively. Let \(f\colon X\to Y\) be a Lipschitz map. Then it holds that
		\begin{equation}\label{eq:area_formula}
			\int_X \eta(x){\rm J}_n^xf\,\d\H^n(x)=\int_Y\sum_{x\in f^{-1}(y)}\eta(x)\,\d\H^n(y)
		\end{equation}
		for every Borel function \(\eta\colon X\to[0,+\infty)\). In particular, we have that \({\rm J}_n^x f\)
		is \(\H^n\)-a.e.\ independent of the chosen atlases \(\mathcal A_X\) and \(\mathcal A_Y\).
	\end{theorem}
	\begin{proof}
		Fix \(i,j\in\N\). By applying the area formula \cite[Theorem 5.1]{Ambrosio-Kirchheim00} to the Lipschitz maps
		\(\varphi^{-1}_i\colon\varphi_i(G_f^{i,j})\to X\) and \(f\circ\varphi_i^{-1}\colon\varphi_i(G_f^{i,j})\to Y\)
		we obtain that
		\[\begin{split}
			\int_{\varphi_i(G_f^{i,j})}\eta_1(\varphi_i^{-1}(z)){\rm J}_n(\|\cdot\|_{\varphi_i^{-1}(z)})\,\d\L^n(z)
			&=\int_{G_f^{i,j}}\eta_1(x)\,\d\H^n(x),\\
			\int_{\varphi_i(G_f^{i,j})}\eta_2(z){\rm J}_n({\rm md}_z(f\circ\varphi_i^{-1}))\,\d\L^n(z)
			&=\int_Y\sum_{z\in\varphi_i(f^{-1}(y)\cap G_f^{i,j})}\eta_2(z)\,\d\H^n(y)
		\end{split}\]
		for every Borel functions \(\eta_1\colon X\to\R\) and \(\eta_2\colon\R^n\to\R\). Using also
		Lemma \ref{lemma:Jacobians}, we get that
		\[\begin{split}
			\int_{G_f^{i,j}}\eta(x){\rm J}_n^x f\,\d\H^n(x)&=\int_{\varphi_i(G_f^{i,j})}\eta(\varphi_i^{-1}(z))
			{\rm J}_n^{\varphi_i^{-1}(z)}f\,{\rm J}_n(\|\cdot\|_{\varphi_i^{-1}(z)})\,\d\L^n(z)\\
			&=\int_{\varphi_i(G_f^{i,j})}\eta(\varphi_i^{-1}(z)){\rm J}_n({\rm md}_z(f\circ\varphi_i^{-1}))\,\d\L^n(z)\\
			&=\int_Y\sum_{z\in\varphi_i(f^{-1}(y)\cap G_f^{i,j})}\eta(\varphi_i^{-1}(z))\,\d\H^n(y)\\
			&=\int_Y\sum_{x\in f^{-1}(y)\cap G_f^{i,j}}\eta(x)\,\d\H^n(y).
		\end{split}\]
		Summing up the above identity over \(i,j\in\N\), we obtain the sought identity \eqref{eq:area_formula}.
	\end{proof}
	\subsection{Porous sets and locality of pointwise Lipschitz constants}\label{sec:localitylip}
	Let \((X,d_X)\) and \((Y,d_Y)\) be metric spaces. Let \(f\colon X\to Y\) be a Lipschitz map.
	We define the \textbf{pointwise Lipschitz constant} function \(\lip(f)\colon X\to[0,+\infty)\) of \(f\) as
	\[
	\lip(f)(x)\coloneqq\limsup_{\tilde x\to x}\frac{d_Y(f(\tilde x),f(x))}{d_X(\tilde x,x)}
	\]
	for every accumulation point \(x\in X\), and \(\lip(f)(x)\coloneqq 0\) for every isolated point \(x\in X\).
	\begin{definition}
		A subset \(P\subseteq X\) of a metric space \((X,d)\) is said to be \textbf{porous} provided
		\[
		\forall x\in P\quad\exists\,\varepsilon>0:\quad\forall r>0\quad\exists\,y\in B_r(x)\setminus\{x\}:\quad d(y,P)\geq\varepsilon\,d(y,x).
		\] 
	\end{definition}
	
	\begin{remark}\label{rmk:equiv_porous_with_dist}{\rm
			Observe that, given any set \(P\subseteq X\), it holds that
			\[
			P\text{ is porous}\quad\Longleftrightarrow\quad\lip(d(\cdot,P))(x)>0\text{ for every }x\in P.
			\]
			Indeed, if \(P\) is porous, \(x\in P\) and \(\varepsilon>0\) is such that for every \(r>0\) there
			exists \(y\in B_r(x)\) with \(d(y,P)\geq\varepsilon\,d(y,x)\), then
			\[
			\lip(d(\cdot,P))(x)=\limsup_{y\to x}\frac{d(y,P)}{d(y,x)}\geq\varepsilon>0.
			\]
			Conversely, if \(\limsup_{y\to x}d(y,P)/d(y,x)=\lip(d(\cdot,P))(x)>0\) for every \(x\in P\), then
			for any \(x\in P\) and \(r>0\) we can find \(y\in B_r(x)\) such that \(d(y,P)\geq\varepsilon_x\,d(y,x)\),
			where we set \(\varepsilon_x\coloneqq\lip(d(\cdot,P))(x)/2>0\), thus proving that \(P\) is porous.
	}\end{remark}
	\begin{lemma}\label{lemma:porous1}
		Let \((X,d_X)\), \((Y,d_Y)\) be metric spaces. Let \(f\colon X\to Y\) be \emph{locally biLipschitz},
		i.e.\ for any \(x\in X\) there exists \(r>0\) such that \(f\colon B_r(x)\to f(B_r(x))\) is biLipschitz. Assume that every porous
		subset of \(Y\) is \(\H^n\)-negligible, for some given \(n\in\N\). Then every porous Borel subset of \(X\) is \(\H^n\)-negligible.
	\end{lemma}
	\begin{proof}
		We can find \((x_k)_{k\in\N}\subset X\), \((r_k)_{k\in\N}\subset(0,1)\), and \((C_k)_{k\in\N}\subset[1,+\infty)\)
		such that \(X=\bigcup_{k\in\N}B_{r_k}(x_k)\) and \(f\colon B_{2 r_k}(x_k)\to f(B_{2 r_k}(x_k))\)
		is \(C_k\)-biLipschitz for every \(k\in\N\). Fix a porous set \(P\subset X\) and denote
		\(P_k\coloneqq P\cap B_{r_k}(x_k)\) for every \(k\in\N\). For any \(k\in\N\), \(x\in P_k\) and
		\(r\in(0,r_k)\) we have \(B_{r/C_k}(x)\subset B_{2r_k}(x_k)\), so that \(f(B_{r/C_k}(x))\subset B_r(f(x))\) and thus
		\[\begin{split}
			{\rm lip}(d_Y(\cdot,f(P_k)))(f(x))&=\inf_{r<r_k}\sup_{y\in B_r(f(x))\setminus\{f(x)\}}\frac{d_Y(y,f(P_k))}{d_Y(y,f(x))}\\
			&\geq\inf_{r<r_k}\sup_{\tilde x\in B_{r/C_k}(x)\setminus\{x\}}\frac{d_Y(f(\tilde x),f(P_k))}{d_Y(f(\tilde x),f(x))}\\
			&\geq\frac{1}{C_k^2}\inf_{r<r_k}\sup_{\tilde x\in B_{r/C_k}(x)\setminus\{x\}}\frac{d_X(\tilde x,P_k)}{d_X(\tilde x,x)}
			=\frac{\lip(d_X(\cdot,P_k))(x)}{C_k^2}
		\end{split}\]
		for every \(x\in P_k\). Taking Remark \ref{rmk:equiv_porous_with_dist} into account, we deduce that
		\(f(P_k)\) is a porous subset of \(Y\). Therefore, we have \(\H^n(P_k)=\H^n(f^{-1}(f(P_k))\cap B_{2r_k}(x_k))
		\leq C_k^n\H^n(f(P_k))=0\), thus \(\H^n(P)\leq\sum_{k\in\N}\H^n(P_k)=0\). This shows that \(\H^n(P)=0\)
		for every porous set \(P\subset X\).
	\end{proof}

	In this paper, by a \textbf{metric measure space} \((X,d,\mathfrak m)\) we mean a metric space
	\((X,d)\) together with a Radon measure \(\mathfrak m\geq 0\) on \(X\) that is finite on balls.
	\begin{definition}
		We say that a metric measure space \((X,d,\mathfrak m)\) is \textbf{porosity-null} if every porous Borel
		set \(P\subseteq X\) satisfies \(\mathfrak m(P)=0\).
	\end{definition}

	Every porosity-null metric measure space \((X,d,\mathfrak m)\) is \textbf{asymptotically doubling}, i.e.
	\[
	\limsup_{r\searrow 0}\frac{\mathfrak m(B_{2r}(x))}{\mathfrak m(B_r(x))}<+\infty\quad\text{ for }\mathfrak m\text{-a.e.\ }x\in X,
	\]
	see \cite[Theorem 3.6]{Mera-Moran-Preiss03}, thus in particular the Lebesgue differentiation theorem holds in \((X,d,\mathfrak m)\). 
	The next example shows that the converse implication does not hold, meaning that an asymptotically doubling metric measure
	space is not necessarily porosity-null.
	\begin{example}{\rm
			Equip \([0,1]^2\) with the Euclidean distance \(d\) and the measure \(\mathcal H^2+\mathcal H^1|_S\), where \(S\) denotes the
			segment \([0,1]\times\{1/2\}\). This space is asymptotically doubling, but \(S\) is a porous set (e.g.\ by Remark
			\ref{rmk:equiv_porous_with_dist}, since \(\lip(d(\cdot,S))=1\) in \(S\)) having positive measure.
	}\end{example}

	\begin{proposition}[Locality of \(\lip\)]\label{prop:equiv_locality_lip}
		Let \((X,d,\mathfrak m)\) be a metric measure space. Then the following conditions are equivalent:
		\begin{itemize}
			\item[\(\rm i)\)] \((X,d,\mathfrak m)\) is porosity-null.
			\item[\(\rm ii)\)] Given any Borel set \(E\subseteq X\), for \(\mathfrak m\)-a.e.\ \(x\in E\) it holds that
			\begin{equation}\label{eq:equiv_porous_null}
				\forall\varepsilon>0\quad\exists\,r>0:\quad\forall y\in B_r(x)\setminus\{x\},\quad B_{\varepsilon\,d(y,x)}(y)\cap E\neq\varnothing.
			\end{equation}
			\item[\(\rm iii)\)] Given any Borel set \(E\subseteq X\) and any \(f\colon X\to Y\) Lipschitz,
			it holds that
			\[
			\lip(f)(x)=\lip(f|_E)(x)\quad\text{ for }\mathfrak m\text{-a.e.\ }x\in E.
			\]
			\item[\(\rm iv)\)] Given two Lipschitz functions \(f,g\colon X\to\R\), it holds that
			\[
			\lip(f)(x)=\lip(g)(x)\quad\text{ for }\mathfrak m\text{-a.e.\ }x\in\{f=g\}.
			\]
		\end{itemize}
	\end{proposition}
	\begin{proof}
		First of all, observe that for any \(E\subseteq X\), \(\varepsilon>0\) and \(x,y\in X\) we have that
		\begin{equation}\label{eq:equiv_locality_lip_aux}
			d(y,E)\geq\varepsilon\,d(y,x)\quad\Longleftrightarrow\quad B_{\varepsilon\,d(y,x)}(y)\cap E=\varnothing.
		\end{equation}
		\(\boldsymbol{{\rm i)}\Longrightarrow{\rm ii)}.}\) Assume ii) does not hold. Then we find Borel sets
		\(P\subseteq E\subseteq X\) with \(\mathfrak m(P)>0\) and \(\varepsilon>0\) such that for any \(r>0\)
		there exists \(y\in B_r(x)\) for which \(B_{\varepsilon\,d(y,x)}(y)\cap E=\varnothing\). In particular
		\(B_{\varepsilon\,d(y,x)}(y)\cap P=\varnothing\), so that \(P\) is porous by \eqref{eq:equiv_locality_lip_aux}.
		All in all, \(P\) is a porous Borel set with \(\mathfrak m(P)>0\), thus accordingly i) does not hold.
		\smallskip
		
		\noindent\(\boldsymbol{{\rm ii)}\Longrightarrow{\rm i)}.}\) Assume ii) holds and let \(P\) be a porous Borel set.
		Applying \eqref{eq:equiv_porous_null} to \(P\) and recalling \eqref{eq:equiv_locality_lip_aux}, we deduce that
		\(\mathfrak m(P)=0\), thus proving i).
		\smallskip
		
		\noindent\(\boldsymbol{{\rm ii)}\Longrightarrow{\rm iii)}.}\) Assume ii) holds. Fix \(E\subseteq X\) Borel and
		\(f\colon X\to Y\) Lipschitz. Clearly, \(\lip(f|_E)(x)\leq\lip(f)(x)\) for all \(x\in E\), thus it suffices to
		show the converse inequality. Take a point \(x\in E\) for which \eqref{eq:equiv_porous_null} holds. If \(x\) is
		an isolated point then \(\lip(f|_E)(x)=\lip(f)(x)\), thus assume \(x\) is an accumulation point. Pick a sequence
		\((x_n)_n\subseteq X\setminus\{x\}\) converging to \(x\) such that \(d(f(x_n),f(x))/d(x_n,x)\to\lip(f)(x)\).
		Applying \eqref{eq:equiv_porous_null} we obtain that, up to a subsequence in \(n\), there exists
		a sequence \((y_n)_n\subseteq E\) such that \(d(x_n,y_n)<d(x_n,x)/n\) for every \(n\in\N\). In particular,
		\(d(y_n,x)\leq d(y_n,x_n)+d(x_n,x)<(1+n^{-1})d(x_n,x)\to 0\) as \(n\to\infty\), and
		\(d(y_n,x)\geq d(x_n,x)-d(y_n,x_n)>(1-n^{-1})d(x_n,x)>0\) for all \(n\geq 2\). Hence,
		\[\begin{split}
			\lip(f)(x)&\leq\limsup_{n\to\infty}\frac{d(f(x_n),f(y_n))}{d(x_n,y_n)}\frac{d(x_n,y_n)}{d(x_n,x)}
			+\limsup_{n\to\infty}\frac{d(f(y_n),f(x))}{d(y_n,x)}\frac{d(y_n,x)}{d(x_n,x)}\\
			&\leq\Lip(f)\lim_{n\to\infty}\frac{1}{n}+\lim_{n\to\infty}\bigg(1+\frac{1}{n}\bigg)\limsup_{n\to\infty}\frac{d(f(y_n),f(x))}{d(y_n,x)}
			\leq\lip(f|_E)(x).
		\end{split}\]
		\(\boldsymbol{{\rm iii)}\Longrightarrow{\rm iv)}.}\) Assume iii) holds and fix \(f,g\colon X\to\R\)
		Lipschitz. Then
		\[
		\lip(f)(x)=\lip(f|_{\{f=g\}})(x)=\lip(g|_{\{f=g\}})(x)=\lip(g)(x)\quad\text{ for }\mathfrak m\text{-a.e.\ }x\in\{f=g\},
		\]
		thus proving iv).
		\smallskip
		
		\noindent\(\boldsymbol{{\rm iv)}\Longrightarrow{\rm i)}.}\) Assume i) does not hold. Then there exists a porous Borel
		set \(P\subseteq X\) such that \(\mathfrak m(P)>0\). Remark \ref{rmk:equiv_porous_with_dist} ensures that
		\(\lip(d(\cdot,P))(x)>0=\lip(0)(x)\) for every \(x\in P\). Since the Lipschitz functions \(d(\cdot,P)\) and \(0\)
		agree in the set \(P\subseteq\{d(\cdot,P)=0\}\), it follows that item iv) does not hold.
	\end{proof}

	We recall that a metric measure space \((X,d,\mathfrak m)\) is said to be \textbf{locally doubling} if for any
	\(\bar x\in X\) there exist constants \(C>0\) and \(\bar r>0\) such that
	\[
	\mathfrak m(B_{2r}(x))\leq C\mathfrak m(B_r(x))\quad\text{ for every }x\in B_{\bar r}(\bar x)\text{ and }r\in(0,\bar r).
	\]
	Every locally-doubling metric measure space satisfies the conditions in Proposition \ref{prop:equiv_locality_lip};
	see e.g.\ \cite[Lemma 1.3]{Gigli-Pasqualetto22}.
	\subsection{Jacobians and pointwise Lipschitz constants}\label{sec:Jacobian_equivalence}
	Throughout this subsection, we assume that:
	\begin{itemize}
		\item \((X,d_X)\) and \((Y,d_Y)\) are \(n\)-rectifiable metric spaces;
		\item \(\H^n_X\) is boundedly finite, i.e.\ \(\H^n(B_r(x))<\infty\) for every \(x\in X\) and \(r>0\);
		\item \(f\colon X\to Y\) is a Lipschitz map;
		\item porous subsets of \(X\) are \(\H^n\)-negligible.
	\end{itemize}
	In particular, \((X,d_X,\H^n)\) is a porosity-null metric measure space, thus in particular
	the Lebesgue differentiation theorem holds in \((X,d_X,\H^n)\).
	\begin{definition}
		Under the above assumptions, we define the \textbf{metric Jacobian} of \(f\) as
		\[
		{\rm Jac}_x f\coloneqq\limsup_{r\searrow 0}\frac{\H^n(f(B_r(x)))}{\H^n(B_r(x))}\quad\text{ for every }x\in X.
		\]
	\end{definition}
	\begin{theorem}\label{thm:Jacobians}
		Under the above assumptions, it holds that
		\begin{equation}\label{eq:Jacobians1}
			{\rm Jac}_x f={\rm J}_n^x f\leq(\lip f)^n(x)\quad\text{ for }\H^n\text{-a.e.\ }x\in X.
		\end{equation}
		Moreover, it holds that
		\begin{equation}\label{eq:Jacobians2}
			\lip f(x)=\|\d_x f\|_{\rm op}\coloneqq\sup_{v\in T_x X\setminus\{0\}}\frac{\|\d_x f(v)\|_{f(x)}}{\|v\|_x}
			\quad\text{ for }\H^n\text{-a.e.\ }x\in X_f,
		\end{equation}
		and \(\|\d_x f\|^n_{\rm op}={\rm J}_n^x f\) if and only if
		\(\frac{\d_x f}{\|\d_x f\|_{\rm op}}\colon T_x X\to T_{f(x)}Y\) is an isometry.
	\end{theorem}
	\begin{proof}
		Let \((A_k)_{k\in\N}\) be as in Lemma \ref{lemma:injectivePartition}. Given any \(x\in X\) and \(r>0\),
		we have that \(\#(f^{-1}(y)\cap B_r(x))\geq 1\) for every \(y\in f(B_r(x))\), thus by applying the Lebesgue
		differentiation theorem and the area formula (Theorem \ref{thm:area_formula}), for \(\H^n\)-a.e.\ \(x\in X\)
		we obtain that
		\[\begin{split}
			{\rm J}_n^x f&=\lim_{r\searrow 0}\fint_{B_r(x)}{\rm J}_n^z f\,\d\H^n(z)
			=\lim_{r\searrow 0}\frac{1}{\H^n(B_r(x))}\int_Y\#(f^{-1}(y)\cap B_r(x))\,\d\H^n(y)\\
			&\geq\limsup_{r\searrow 0}\frac{\H^n(f(B_r(x)))}{\H^n(B_r(x))}={\rm Jac}_x f.
		\end{split}\]
		Conversely, given any \(k\in\N\), \(x\in A_k\) and \(r>0\), we have that \(\#(f^{-1}(y)\cap B_r(x)\cap A_k)=1\)
		for every \(y\in f(B_r(x)\cap A_k)\), thus a similar computation to the one above yields
		\[\begin{split}
			{\rm J}_n^x f&=\lim_{r\searrow 0}\fint_{B_r(x)\cap A_k}{\rm J}_n^z f\,\d\H^n(z)
			=\lim_{r\searrow 0}\frac{1}{\H^n(B_r(x))}\int_Y\#(f^{-1}(y)\cap B_r(x)\cap A_k)\,\d\H^n(y)\\
			&=\lim_{r\searrow 0}\frac{\H^n(f(B_r(x)\cap A_k))}{\H^n(B_r(x))}
			\leq\limsup_{r\searrow 0}\frac{\H^n(f(B_r(x)))}{\H^n(B_r(x))}={\rm Jac}_x f.
		\end{split}\]
		for \(\H^n\)-a.e.\ density point \(x\) of \(A_k\). All in all, we have that
		\({\rm J}_n^x f={\rm Jac}_x f\) for \(\H^n\)-a.e.\ \(x\in X\).
		
		Fix \(i,j\in\N\). For \(\H^n\)-a.e.\ \(x\in X_f^{i,j}\), Proposition \ref{prop:equiv_locality_lip}
		and the relevant definitions yield
		\[\begin{split}
			\lip(f)(x)&=\lip(f|_{X_f^{i,j}})(x)=\limsup_{X_f^{i,j}\ni z\to x}\frac{d_Y(f(z),f(x))}{d_X(z,x)}
			=\limsup_{X_f^{i,j}\ni z\to x}\frac{\|\psi_j(f(z))-\psi_j(f(x))\|_{f(x)}}{d_X(z,x)}\\
			&=\limsup_{X_f^{i,j}\ni z\to x}\frac{\|\d_x f(\varphi_i(z)-\varphi_i(x))\|_{f(x)}}{d_X(z,x)}\\
			&=\limsup_{X_f^{i,j}\ni z\to x}\bigg\|\d_x f\bigg(\frac{\varphi_i(z)-\varphi_i(x)}
			{\|\varphi_i(z)-\varphi_i(x)\|_x}\bigg)\bigg\|_{f(x)}\frac{\|\varphi_i(z)-\varphi_i(x)\|_x}{d_X(z,x)}\\
			&=\limsup_{X_f^{i,j}\ni z\to x}\bigg\|\d_x f\bigg(\frac{\varphi_i(z)-\varphi_i(x)}
			{\|\varphi_i(z)-\varphi_i(x)\|_x}\bigg)\bigg\|_{f(x)}.
		\end{split}\]
		Whenever \(\varphi_i(x)\) is a density point of \(\varphi_i(X_f^{i,j})\), we have that every unit vector
		in \(T_x X\) can be approximated by vectors of the form
		\(\frac{\varphi_i(z)-\varphi_i(x)}{\|\varphi_i(z)-\varphi_i(x)\|_x}\) with \(z\in X_f^{i,j}\). Therefore,
		the above computations show that \(\lip(f)(x)=\|\d_x f\|_{\rm op}\) for \(\H^n\)-a.e.\ \(x\in X_f\),
		thus proving \eqref{eq:Jacobians2}.
		
		Now, observe that for \(\H^n\)-a.e.\ \(x\in X_f\) it follows from the very definition of \(\|\d_x f\|_{\rm op}\) that
		\begin{equation}\label{ineq:proofprop1}
			\big\{v\in T_x X\;\big|\;\|\d_x f(v)\|_{f(x)}\leq 1\big\}\supset
			\big\{v\in T_x X\;\big|\;\|v\|_x\leq \|\d_x f\|_{\rm op}^{-1}\big\}.
		\end{equation}
		Since \cite[Lemma 6(i)]{Kirchheim94} ensures that \(\H^n(\{v\in T_x X:\|v\|_x\leq 1\})=\omega_n\), we deduce that
		\[\begin{split}
			{\rm J}_n^x f&=\frac{\omega_n}{\H^n(\{v\in T_x X:\|\d_x f(v)\|_{f(x)}\leq 1\})}
			\leq\frac{\omega_n}{\H^n(\{v\in T_x X:\|v\|_x\leq\|\d_x f\|_{\rm op}^{-1}\})}\\
			&=\frac{\omega_n}{\|\d_x f\|_{\rm op}^{-n}\H^n(\{v\in T_x X:\|v\|_x\leq 1\})}=\|\d_x f\|_{\rm op}^n=(\lip f)^n(x)
		\end{split}\]
		for \(\H^n\)-a.e.\ \(x\in X_f\), thus completing the proof of \eqref{eq:Jacobians1}. Furthermore, it follows from
		the above estimates that the equality \({\rm J}_n^x f=\|\d_x f\|_{\rm op}^n\) is achieved if and only if
		\begin{equation}\label{ineq:proofprop2}
			\H^n\big(\big\{v\in T_x X\;\big|\;\|\d_x f(v)\|_{f(x)}\leq 1\big\}\setminus
			\big\{v\in T_x X\;\big|\;\|v\|_x\leq \|\d_x f\|_{\rm op}^{-1}\big\}\big)=0.
		\end{equation}
		In view of the continuity of \(v\mapsto\|\d_x f(v)\|_{f(x)}\) and \(v\mapsto\|v\|_x\), we have that
		\eqref{ineq:proofprop2} is equivalent to 
		\(\big\{v\in T_x X\;\big|\;\|\d_x f(v)\|_{f(x)}\leq 1\big\}=\big\{v\in T_x X\;\big|\;\|v\|_x\leq \|\d_x f\|_{\rm op}^{-1}\big\}\), which in turn is equivalent to \(\d_x f(B_{\|\d_x f\|_{\rm op}^{-1}}(0))=B_1(0)\).
		By scaling, the latter can be rewritten as
		\[
		\d_x f(B_1(0))=B_{\|\d_x f\|_{\rm op}}(0).
		\]
		The final claim about on the isometricity of \(\frac{\d_x f}{\|\d_x f\|_{\rm op}}\) then follows
		from the fact that both this normalized operator and its inverse have operator norm equal to \(1\).
	\end{proof}

	For any Lipschitz function \(\phi\colon X\to\R\), by arguing as in the proof of \cref{thm:Jacobians} one can prove
	the following statement.
	\begin{corollary}\label{lemma:dualnorm-lip}
		Under the above assumptions, it holds that
		\[
		\lip(\phi)(x)=\|\tilde\d_x\phi\|_x^*\coloneqq\sup_{w\in T_x X\setminus\{0\}}\frac{\tilde\d_x\phi(w)}{\|w\|_x}
		\quad\text{ for }\H^n\text{-a.e.\ }x\in X.
		\]
	\end{corollary}
	\subsection{Volume-preserving 1-Lipschitz maps between rectifiable spaces}\label{sec:LVR_rectifiable1}
	Let us begin with some auxiliary results that will be used to prove our first main result.
	\begin{lemma}\label{lemma:volume-preserving}
		Let \((X,d_X)\), \((Y,d_Y)\) be metric spaces. Assume that \(\H^n_X\) is \(\sigma\)-finite,
		i.e.\ there exist Borel subsets \(\{E_i:i\in\N\}\) of \(X\) such that \(\H^n_X(E_i)<+\infty\) for
		every \(i\in\N\) and \(X=\bigcup_{i\in\N}E_i\). Let \(f\colon X\to Y\) be \(1\)-Lipschitz. Then \(f(A)\)
		is \(\H_Y^n\)-measurable for every Borel set \(A\subset X\). Moreover, the following conditions are equivalent:
		\begin{itemize}
			\item[\(\rm i)\)] \(f_\#\mathcal H^n_X=\mathcal H^n_Y|_{f(X)}\), which means that
			\[
			\H^n_Y(B\cap f(X))=\H^n_X(f^{-1}(B))\quad\text{ for every Borel set }B\subset Y;
			\]
			\item[\(\rm ii)\)] $\H_Y^n(f(A))=\H_X^n(A)$ for every Borel set $A\subset X$.
		\end{itemize}
		If in addition $\H_X^n(X)<\infty$, then these conditions are equivalent to $\H_Y^n(f(X))=\H_X^n(X)$.
	\end{lemma}
	\begin{proof}
		Since \(\H^n_X\) is \(\sigma\)-finite, every Borel subset \(A\) of
		\(X\) is \(\H^n_X\)-equivalent to some \(\sigma\)-compact set \(\tilde A\subset A\). Since \(f\) is \(1\)-Lipschitz, we get that
		\[
		\H^n_Y(f(A)\setminus f(\tilde A))\leq\H^n_Y(f(A\setminus\tilde A))\leq\H^n_X(A\setminus\tilde A)=0,
		\]
		so that \(f(A)\) is \(\H^n_Y\)-equivalent to the \(\sigma\)-compact set \(f(\tilde A)\). Hence, \(f(A)\) is \(\H_Y^n\)-measurable.
		
		To prove that i) implies ii), assume that $f_\#\H^n_X=\H^n_Y|_{f(X)}$. Then we have that
		\[
		\H_X^n(A)\leq\H_X^n(f^{-1}(f(A)))=\H_Y^n(f(A)\cap f(X))=\H^n_Y(f(A))\leq\H^n_X(A),
		\]
		where we used the inclusion \(A\subset f^{-1}(f(A))\) and the \(1\)-Lipschizianity of the map \(f\).
		
		To prove that ii) implies i), assume that $\H_Y^n(f(A))=\H_X^n(A)$ for every Borel set $A\subset X$.
		For any Borel set \(B\subset Y\) we have that \(f(f^{-1}(B))=B\cap f(X)\), thus accordingly
		\[
		\H^n_Y(B\cap f(X))=\H^n_Y(f(f^{-1}(B)))=\H^n_X(f^{-1}(B)).
		\]
		Finally, if \(\H^n_X(X)<\infty\), then \(\H^n_Y(f(X))=\H^n_X(X)\) implies for every Borel set \(A\subset X\) that
		\[\begin{split}
			\H^n_X(X)&=\H^n_Y(f(X))=\H^n_Y(f(A)\cup f(X\setminus A))\leq\H^n_Y(f(A))+\H^n_Y(f(X\setminus A))\\
			&\leq\H^n_X(A)+\H^n_X(X\setminus A)=\H^n_X(X),
		\end{split}\]
		which forces the identities \(\H^n_Y(f(A))=\H^n_X(A)\) and \(\H^n_Y(f(X\setminus A))=\H^n_X(X\setminus A)\).
	\end{proof}
	
	\begin{proposition}\label{prop:LVR}
		Let \((X,d_X)\), \((Y,d_Y)\) be \(n\)-rectifiable metric spaces.
		Assume that \(\H^n_X\), \(\H^n_Y\) are boundedly finite and \((X,d_X,\H^n_X)\), \((Y,d_Y,\H^n_Y)\) are porosity-null.
		Let \(f\colon X\to Y\) be a \(1\)-Lipschitz map such that \(f_\#\H^n_X=\H^n_Y|_{f(X)}\). Then:
		\begin{enumerate}
			\item \textsc{Essential injectivity.} \(f^{-1}(f(x))=\{x\}\) for \(\H^n\)-a.e.\ \(x\in X\);
			\item \textsc{Chain rule.} \(\lip(\phi\circ f)=\lip(\phi)\circ f\) holds \(\H^n\)-a.e.\ on \(X\) for every \(\phi\in\Lip(Y)\).
		\end{enumerate}
	\end{proposition}
	\begin{proof}
		Note that \(\H^n(f(A))=\H^n(A)\) for every Borel set \(A\subset X\) by \cref{lemma:volume-preserving}, thus in particular
		\({\rm Jac}_x f\geq 1\) for every \(x\in X\). Since \(\lip f\leq 1\) by the \(1\)-Lipschitzianity of \(f\),
		we deduce from \cref{thm:Jacobians} that
		\begin{equation}\label{eq:proofThm_lip2}
			{\rm Jac}_x f={\rm J}_n^x f=(\lip f)^n(x)=1\quad\text{ for }\H^n\text{-a.e.\ }x\in X
		\end{equation}
		and the tangential differential \(\d_x f\colon T_x X\to T_{f(x)}Y\) is an isometry for \(\H^n\)-a.e.\ \(x\in X_f\).
		It also follows from \eqref{eq:proofThm_lip2} that \(\H^n(X\setminus X_f)=0\). Taking also \cref{lemma:dualnorm-lip},
		\cref{lemma:chainrule} and the fact that \(f_\#\H^n_X\ll\H^n_Y\) into account, we then deduce that
		\[
		\lip(\phi\circ f)(x)=\|\tilde\d_x(\phi\circ f)\|^*_x=\|\tilde\d_{f(x)}\phi\circ\d_x f\|_x^*
		=\|\tilde\d_{f(x)}\phi\|_{f(x)}^*=\lip(\phi)(f(x))
		\]
		for \(\H^n\)-a.e.\ \(x\in X\), thus proving (2). To prove (1), it suffices to show that for any Borel set
		\(E\subset X\) with \(\H^n(E)<\infty\) it holds that \(\#(f^{-1}(f(x))\cap E)=1\) for \(\H^n\)-a.e.\ \(x\in E\). Applying \eqref{eq:proofThm_lip2} and \cref{thm:area_formula}, we obtain that
		\[\begin{split}
			\H^n(E)&=\int_E{\rm J}_n^x f\,\d\H^n(x)=\int_Y\#(f^{-1}(y)\cap E)\,\d\H^n(y)
			=\int_{f(E)}\#(f^{-1}(y)\cap E)\,\d\H^n(y)\\
			&\geq\H^n(f(E))=\H^n(E).
		\end{split}\]
		Therefore, the inequality in the above formula is an equality, so that \(\#(f^{-1}(y)\cap E)=1\) for
		\(\H^n_Y\)-a.e.\ \(y\in f(E)\), and accordingly \(\#(f^{-1}(f(x))\cap E)=1\) for \(\H^n\)-a.e.\ \(x\in E\).
	\end{proof}

	We now introduce another condition on metric measure spaces that will be useful in what follows.
	We defer a more detailed discussion of this condition to Section \ref{app:lip-to-Lip}.
	\begin{definition}\label{def:L}
		We say that a metric measure space \((X,d,\mathfrak m)\) has the \textbf{lip-to-Lip property} provided it satisfies
		\[
		{\rm Lip}(\phi)=\|{\rm lip}(\phi)\|_{L^\infty(\mathfrak m)}\quad \text{ for every }\phi\in\Lip(X).
		\]
	\end{definition}
	\medskip
	Given two metric spaces \((X,d_X)\) and \((Y,d_Y)\), we say that a map \(f\colon X\to Y\) is
	\begin{itemize}
		\item a \textbf{local isometry} if for any $x\in X$ there exists an open neighborhood $U$ of $x$ such that $f|_U\colon U\to f(U)$ is an isometry,
		\item a \textbf{length isometry} if for any rectifiable curve $\gamma$ in $X$ we have $L_X(\gamma)=L_Y(f\circ \gamma)$.
	\end{itemize}
	Here, the length of a continuous curve $\gamma\colon [0,1]\to X$ is defined by
	\begin{equation}
		L_X(\gamma)\coloneqq\sup\bigg\{\sum_{i=1}^k d_X(\gamma_{t_{i-1}},\gamma_{t_i})\;\bigg|\;0=t_0<\ldots<t_k=1\bigg\},
	\end{equation}
	where the supremum is taken over all finite partitions \(t_0,\ldots,t_k\) of $[0,1]$, and by a rectifiable curve we mean
	a continuous curve with finite length.
	\begin{theorem}\label{thm:LVR-rectifiable}
		Assume that:
		\begin{enumerate}[label=(\arabic*)]
			\item \((X,d_X)\) is an \(n\)-rectifiable metric space with \(\H^n_X\) boundedly finite.
			\item \((Y,d_Y)\) is a metric space with \(\H^n_Y\) boundedly finite such that
			\((Y,d_Y,\H^n_Y)\) is porosity-null and has the lip-to-Lip property.
			\item\label{item:LVR3} $f\colon X\to Y$ is an open, locally biLipschitz, $1$-Lipschitz map satisfying \(f_\#\H_X^n=\H_Y^n|_{f(X)}.\)
		\end{enumerate}
		Then \(f\) is a local isometry and a length isometry. Moreover, if in addition \((X,d_X)\) and \((Y,d_Y)\)
		are complete, then \(f\colon X\to Y\) is an isometric bijection.
	\end{theorem}
	\begin{proof}
		Fix any \(x_0\in X\). By assumption, there exists an open set \(U\subset X\) with \(x_0\in U\) such that \(V\coloneqq f(U)\) is open
		in \(Y\), and \(f\colon U\to V\) is \(C\)-biLipschitz for some constant \(C\geq 1\). Since \(U\) is \(n\)-rectifiable
		and \(V\) is biLipschitz equivalent to \(U\), we have that \(V\) is \(n\)-rectifiable.
		Lemma \ref{lemma:porous1} ensures that \((U,d_X,\H^n)\) is porosity-null, whereas Lemma \ref{lemma:volume-preserving} implies that
		\(f_\#\H^n_U=\H^n_V\). Now, denote \(y_0\coloneqq f(x_0)\in V\) and define \(\delta\coloneqq d_Y(y_0,Y\setminus V)\wedge 1>0\).
		Pick \(0<r<\frac{\delta}{2C+1}\) and call \(U_r\coloneqq f^{-1}(B_r(y_0))\cap U\). Fix \(x_1\in U_r\) and define \(\phi\colon U\to[0,1]\) as
		\[
		\phi(x)\coloneqq d_X(x,x_1)\wedge d_Y(f(x),Y\setminus V)\wedge 1\quad\text{ for every }x\in U.
		\]
		For any \(x\in U_r\), we can estimate \(d_X(x,x_1)\leq Cd_Y(f(x),f(x_1))\leq 2Cr<\delta-r<1\) and
		\[
		d_Y(f(x),Y\setminus V)\geq d_Y(y_0,Y\setminus V)-d_Y(y_0,f(x))>\delta-r.
		\]
		Consequently, we have that \(\phi(x)=d_X(x,x_1)\) for every \(x\in U_r\). Let us define
		\[
		\Phi\coloneqq\phi\circ(f|_U)^{-1}\colon V\to[0,1].
		\]
		Since \(\Phi\) is Lipschitz and \(0\leq\Phi(y)\leq d_Y(y,Y\setminus V)\) for every \(y\in V\), the zero extension
		of \(\Phi\), which we still denote by \(\Phi\colon Y\to[0,1]\), is a Lipschitz function. Applying \cref{prop:LVR}~(2)
		to \(f\colon U\to V\) and exploiting the \(1\)-Lipschitzianity of \(\phi\), we obtain that
		\[
		\lip(\Phi)\circ f=\lip(\Phi\circ f)=\lip\,\phi\leq 1\quad\text{ holds }\H_X^n\text{-a.e.\ in }U,
		\]
		where we do not distinguish between the pointwise Lipschitz constants in \(U\) (resp.\ \(V\)) and in \(X\) (resp.\ \(Y\)),
		as \(U\) and \(V\) are open. Since \(f_\#\H^n_U=\H^n_V\), it follows that \(\lip\,\Phi\leq 1\) holds \(\H_Y^n\)-a.e.\ in \(V\).
		Moreover, if \(y\in Y\setminus V\) and \(z\in Y\) are given, then we have that
		\[
		|\Phi(z)-\Phi(y)|=\Phi(z)\leq d_Y(z,Y\setminus V)\leq d_Y(z,y),
		\]
		whence it follows that \(\lip\,\Phi(y)\leq 1\). All in all, we have shown that \(\lip\,\Phi\leq 1\) holds \(\H^n_Y\)-a.e.\ in \(Y\).
		By the lip-to-Lip property of \((Y,d_Y,\H^n_Y)\), we deduce that \(\Phi\) is \(1\)-Lipschitz on \(Y\). Therefore, for any \(x_2\in U_r\)
		we have that
		\[
		d_X(x_1,x_2)=|\phi(x_1)-\phi(x_2)|=|\Phi(f(x_1))-\Phi(f(x_2))|\leq d_Y(f(x_1),f(x_2))\leq d_X(x_1,x_2).
		\]
		By the arbitrariness of \(x_1,x_2\in U_r\), we have proven that \(f|_{U_r}\colon U_r\to Y\) is an isometry, thus
		\(f\) is a local isometry. Next, we check that \(f\) is a length isometry. To this aim, fix a rectifiable curve
		\(\gamma\colon [0,1]\to X\). Given any \(\varepsilon>0\), we can find \(0=t_0<t_1<\ldots<t_k=1\) such that
		\(d_X(\gamma_{t_{i-1}},\gamma_{t_i})=d_Y(f(\gamma_{t_{i-1}}),f(\gamma_{t_i}))\) for all \(i=1,\ldots,k\) and
		\(\sum_{i=1}^k d_X(\gamma_{t_{i-1}},\gamma_{t_i})>L_X(\gamma)-\varepsilon\). Taking also the \(1\)-Lipschitzianity
		of \(f\) into account, we infer that
		\[
		L_X(\gamma)-\varepsilon<\sum_{i=1}^k d_X(\gamma_{t_{i-1}},\gamma_{t_i})=\sum_{i=1}^k d_Y(f(\gamma_{t_{i-1}}),f(\gamma_{t_i}))\leq L_Y(f\circ\gamma)\leq L_X(\gamma).
		\]
		Letting \(\varepsilon\to 0\), we conclude that \(L_Y(f\circ\gamma)=L_X(\gamma)\), showing that \(f\) is a length isometry.
		\medskip
		
		Let us now assume in addition that \((X,d_X)\) and \((Y,d_Y)\) are complete. By Proposition \ref{prop:(L)-1}, \(Y\) is a length space.
		We denote by \(\mathcal G\) the collection of all open subsets \(U\) of \(X\) such that \(f(U)\subset Y\) is open and \(f|_U:U\to f(U)\) is biLipschitz.
		Since \(f\) is open, \(f(X)\) is open in \(Y\).
		We also claim that \(f\) is injective. To prove it, we argue by contradiction: assume that \(f(x_1)=f(x_2)\) for some \(x_1,x_2\in X\)
		with \(x_1\neq x_2\). Then we can find \(U_1,U_2\in\mathcal G\) such that \(x_1\in U_1\), \(x_2\in U_2\), \(U_1\cap U_2=\varnothing\)
		and \(f(U_1)=f(U_2)\). Since \(f|_{U_2}\colon U_2\to f(U_1)\) is surjective, we have that \(\# f^{-1}(f(x))\geq 2\) for every \(x\in U_1\).
		This contradicts \cref{prop:LVR} (1), as \(\H^n_X(U_1)>0\) due to the fact that \(\H^n_X\) is fully supported.
		Therefore, we have shown that \(f\colon X\to f(X)\) is bijective. Next, we claim that \(f(X)=Y\). Let us argue by contradiction:
		assume there exists \(z\in Y\setminus f(X)\). Fix any \(y\in f(X)\) and take a constant-speed rectifiable curve \(\sigma\colon[0,1]\to Y\)
		such that \(\sigma_0=y\) and \(\sigma_1=z\). Define
		\[
		\delta\coloneqq\inf O,\quad\text{ where we set }O\coloneqq\{t\in[0,1]:\sigma_t\notin f(X)\}.
		\]
		Since \(f(X)\) is open in \(Y\), we have that \(\delta>0\). Note that \(\sigma([0,\delta))\subseteq f(X)\), thus it makes sense to
		define \(\gamma\colon[0,\delta)\to X\) as \(\gamma_t\coloneqq f^{-1}(\sigma_t)\) for every \(t\in[0,\delta)\). Now, fix \(s\in(0,\delta)\).
		The set \(\sigma([0,s])\) is a compact subset of \(f(X)\), thus the Lebesgue number lemma guarantees that we can find
		\(0=t_0<t_1<\ldots<t_k=s\) and \(U_1,\ldots,U_k\in\mathcal G\) such that \(\sigma([t_{i-1},t_i])\subset f(U_i)\) for every \(i=1,\ldots,k\).
		It follows that \(\gamma|_{[t_{i-1},t_i]}\) is a rectifiable curve for every \(i=1,\ldots,k\) and \(\gamma|_{[0,s]}\) is continuous.
		Due to the fact that \(f\) is a length isometry, we deduce that
		\[
		L_X(\gamma|_{[0,s]})=\sum_{i=1}^k L_X(\gamma|_{[t_{i-1},t_i]})=\sum_{i=1}^k L_Y(\sigma|_{[t_{i-1},t_i]})=L_Y(\sigma|_{[0,s]})=sL_Y(\sigma).
		\]
		Therefore, for any \(s,\tilde s\in(0,\delta)\) with \(s<\tilde s\) we can estimate
		\begin{equation}\label{eq:LVR-rectifiable_aux}
			d_X(\gamma_s,\gamma_{\tilde s})\leq L_X(\gamma|_{[s,\tilde s]})=L_X(\gamma|_{[0,\tilde s]})-L_X(\gamma|_{[0,s]})=(\tilde s-s)L_Y(\sigma),
		\end{equation}
		which implies that \(x=\lim_{s\nearrow\delta}\gamma_s\in X\) exists by the completeness of \(X\). We deduce that
		\(\sigma_\delta=\lim_{s\nearrow\delta}\sigma_s=\lim_{s\nearrow\delta}f(\gamma_s)=f(x)\in f(X)\), so that \(\delta\notin O\). We distinguish two cases:
		\begin{itemize}
			\item \(0<\delta<1\). Then \(\sigma_\delta\in f(X)\) and the openness of \(f(X)\) yield \(\sigma([0,\delta+\eta))\subset f(X)\)
			for some \(\eta>0\), or in other words \([0,\delta+\eta)\cap O=\varnothing\), in contradiction with \(\delta=\inf O\).
			\item \(\delta=1\). Then \(z=\sigma_1\in f(X)\), in contradiction with the assumption that \(z\notin f(X)\).
		\end{itemize}
		All in all, we have shown that \(f(X)=Y\), and thus \(f\colon X\to Y\) is a bijection, as we claimed. Finally, let us check that
		\(f\colon X\to Y\) is an isometry. To this aim, fix \(x,\tilde x\in X\) and \(\varepsilon>0\). Since \(Y\) is a length space,
		there exists
		a constant-speed curve \(\sigma\colon[0,1]\to Y\) such that \(\sigma_0=f(x)\), \(\sigma_1=f(\tilde x)\) and \(L_Y(\sigma)\leq d_Y(f(x),f(\tilde x))+\varepsilon\).
		Then the above arguments show that \(\gamma\coloneqq f^{-1}\circ\sigma\) is a rectifiable curve, and the computations in \eqref{eq:LVR-rectifiable_aux} ensure that
		\[
		d_Y(f(x),f(\tilde x))\leq d_X(x,\tilde x)=d_X(\gamma_0,\gamma_1)\leq L_Y(\sigma)\leq d_Y(f(x),f(\tilde x))+\varepsilon.
		\]
		Letting \(\varepsilon\searrow 0\), we conclude that \(d_X(x,\tilde x)=d_Y(f(x),f(\tilde x))\). Hence, \(f\) is an isometry. 
	\end{proof}

	The assumptions in \cref{thm:LVR-rectifiable} may not be optimal.
	Nevertheless, the following counterexamples show that several of the hypotheses cannot be omitted.
	Whether the porosity-null assumption can be removed remains open.
	Below, we collect a list of counterexamples illustrating this, where all measures are understood as the corresponding Hausdorff measures.
	\begin{enumerate}[label=(\alph*), font=\normalfont]
		\item \emph{Local biLipschitzianity.}
		Let \(n\ge2\), and let \(C\subset[0,1]\) be the Cantor set. For each component \((a_j,b_j)\)
		of \([0,1]\setminus C\), let \(D_j\subset\mathbb R^n\) be the closed ball with diameter \([a_je_1,b_je_1]\). Let \(c: [0,1]\to [0,1]\) be the standard Cantor function. Recall that \(c\) is continuous and nondecreasing, is constant on each \([a_j,b_j]\), and satisfies \(c(0)=0\) and
		\(c(3^{-k})=2^{-k}\) for every integer \(k\ge0\). Set
		\[
		\qquad S=\{te_1:0\le t\le 1\},\quad X=Y=S\cup\bigcup_jD_j,
		\quad d_X(p,q)=d_Y(p,q)+|c(p_1)-c(q_1)|,
		\]
		where \(d_Y\) is the intrinsic Euclidean length distance on \(S\cup\bigcup_jD_j\).
		
		The identity \(f=\mathrm{id}:(X,d_X)\to(Y,d_Y)\) is a \(1\)-Lipschitz
		homeomorphism between compact spaces. On each \(D_j\), both distances agree
		with the Euclidean distance. Furthermore, \(S\) is \(\mathcal H^n\)-null for both metrics, which proves rectifiability
		and volume preservation. Since \(Y\) is locally Euclidean outside \(S\cup\bigcup_j\partial D_j\), a null set, it is porosity-null.
		
		For \(\varphi\in\operatorname{Lip}(Y)\), put
		\(L=\|\operatorname{lip}_Y\varphi\|_{L^\infty(\mathcal H_Y^n)}\).
		Each \(\varphi|_{D_j}\) is \(L\)-Lipschitz. The Lipschitz function \(v(t)=\varphi(te_1)\)
		therefore satisfies \(|v'|\le L\) almost everywhere outside \(C\).
		Since \(C\) is Lebesgue-null, absolute continuity makes \(v\)
		\(L\)-Lipschitz on \([0,1]\). Shortest paths between distinct balls
		consist of segments in the endpoint balls and a segment in \(S\).
		Thus \(\operatorname{Lip}(\varphi)\le L\), proving the desired property.
		
		However, \(p_k=3^{-k}e_1\to0\) in \(X\), and
		\[
		\frac{d_X(0,p_k)}{d_Y(0,p_k)}
		=1+\left(\frac32\right)^k\longrightarrow\infty.
		\]
		Thus \(f\) is neither locally biLipschitz nor a local isometry at \(0\).
		
		\item \emph{The lip-to-Lip property}.
		See \cref{example:shortcut-curve} for $\lambda\in (0,1)$.
		The map $f$ is globally biLipschitz, open, but not a local isometry.
		Because the co-domain $Y$ does not satisfy the lip-to-Lip property (actually it is easier to see that $Y$ is not an essentially length space; see \cref{app:lip-to-Lip}).
		\item \emph{Openness}. Let
		\[
		A\coloneqq\R^2\setminus(0,\infty)^2,\qquad B\coloneqq[0,\infty)^2,
		\]
		where $A$ is endowed with the intrinsic metric and $B$ with
		the Euclidean metric. 
		Consider $X\coloneqq A\sqcup B$ the disjoint union, endowed with the corresponding metric on each component and
		\[
		d_X(x,y)\coloneqq |x|+1+|y|,\qquad x\in A,\ y\in B
		\]
		between the two components.
		Let $Y=\R^2$ and let $f\colon X\to Y$ be the natural inclusion on each component.
		
		Then $f$ is volume preserving and $1$-Lipschitz.
		Moreover, it is locally biLipschitz around each $x\in X$ but there is no neighborhood $U$ of $(0,0)\in A$ such that $f|_U$ is biLipschitz and $f(U)$ is open in $Y$.
		\item \emph{Completeness}.
		Consider $f\colon (0,2\pi)\to \mathbb S^1$ by sending $x\mapsto e^{ix}$.
		All conditions are satisfied except for the completeness of $(0,2\pi)$.
		Clearly, $f$ is not a global isometry.
	\end{enumerate}
	
	
	\section{Lipschitz volume rigidity for non-collapsed RCD spaces}
	
	In this section, we establish Lipschitz volume rigidity for non-collapsed RCD spaces
	as stated in Theorem~\ref{thm:main_2}.
	The case $N=1$ follows from the classification of one-dimensional RCD spaces
	\cite[Corollary~1.2]{Kitabeppu-Lakzian16}, with the boundary given by the endpoints.
	We therefore assume $N\geq2$ throughout this section.
	
	\subsection{Preliminaries on RCD spaces}
	A metric measure space $(X,d,\m)$ is said to satisfy the $\rcd(K,N)$ condition, for $K\in\R$ and
	$N\in[1,\infty)$, if it is complete, satisfies the curvature-dimension condition $\mathrm{CD}(K,N)$ in the
	sense of Lott--Sturm--Villani \cite{LottVillani09,Sturm06-2}, and is infinitesimally Hilbertian, namely, the
	Sobolev space $W^{1,2}(X,d,\m)$ is a Hilbert space. We refer e.g.\ to
	\cite{AGS14-RCD,Gigli15-RCD,Gigli-Pasqualetto-Book} for more background on RCD spaces.
	
	Following \cite{DePhilippis-Gigli18}, we call an $\rcd(K,N)$ space $(X,d,\m)$ \textbf{non-collapsed} if
	$\m=\H^N$. In this case, $N$ is necessarily an integer, and by \cite{Mondino-NaberJEMS,DePhilippisMarcheseRindler2017,KellMondino2018,GigliPasqualetto2021} $(X,d,\H^N)$ is
	$N$-rectifiable in the sense of \cref{sec:pre_rectifiable}.
	In particular, all the results from the previous section apply. Unless otherwise stated, throughout this section we additionally
	assume that every non-collapsed RCD space under consideration has full support, i.e., $\supp(\H^N)=X$. We write
	$\mathrm{ncRCD}(K,N)$ for the class of non-collapsed $\rcd(K,N)$ spaces.
	\subsubsection*{Regular/Singular set}
	Let $(X,d,\H^N)$ be a $\mathrm{ncRCD}(K,N)$ space. For $\varepsilon, r > 0$, define the $(\varepsilon,
	r)$-regular set of $X$ by
	\begin{equation}
		\mathcal R_{\varepsilon, r}(X) \coloneqq \{x \in X : \exists t > r \text{ such that }
		d_{\mathrm{GH}}(B_s^X(x), B_s^{\R^N}(0)) \leq \varepsilon s, \text{ for all } s \in (0, t]\}.
	\end{equation}
	The $\varepsilon$-regular set and the regular set of $X$ are then defined respectively by
	\[
	\mathcal R_{\varepsilon}(X) \coloneqq \cup_{r>0}\mathcal R_{\varepsilon, r}(X)\,, \qquad \mathcal R(X)\coloneqq \cap_{\varepsilon>0}\mathcal R_{\varepsilon}(X)\,.
	\]
	By the Bishop--Gromov inequality, the limit
	\[
	\theta_X(x)\coloneqq\lim_{r\searrow 0}\frac{\H^N(B_r(x))}{\omega_N r^N}
	\]
	exists for every $x\in X$ and is called the \textbf{Bishop--Gromov density} of $X$ at $x$. Moreover,
	$\theta_X$ is lower semi-continuous on $X$; see \cite[Lemma 2.2]{DePhilippis-Gigli18}.
	
	By \cite{DePhilippis-Gigli18}, $\theta_X(x)\leq 1$ for every $x\in X$, and $\theta_X(x)=1$ if and only if
	$x\in\mathcal R(X)$. On the other hand, the local Ahlfors \(N\)-regularity of \((X,d,\H^N)\) yields
	\(\theta_X(x)>0\) for every \(x\in X\). In summary,
	\begin{equation}
		\theta_X(x)=
		\begin{cases}
			1, & x\in\mathcal R(X),\\
			\text{a number in }(0,1), & x\in X\setminus\mathcal R(X).
		\end{cases}
	\end{equation}
	
	The continuity of $\H^N$ under Gromov--Hausdorff convergence and the volume rigidity theorem of
	\cite{DePhilippis-Gigli18} yield a quantitative comparison between the deviation of the volume ratio from $1$
	and the Gromov--Hausdorff distance from a Euclidean ball. The corresponding infinitesimal statement is
	recorded in \cite[Corollary 2.14]{Kapovitch-Mondino21}.
	
	Following the convention in \cite{Cheeger-Colding97}, we denote by $\R^j_+ \times \R^k \ni
	(\boldsymbol{\varepsilon}, \mathbf{a}) \mapsto \Psi(\boldsymbol{\varepsilon} \mid \mathbf{a}) \in \R_+$ a
	function s.t.\ for any fixed $\mathbf{a}\in \R^k$, we have $\lim_{\boldsymbol{\varepsilon} \to 0}
	\Psi(\boldsymbol{\varepsilon}\mid  \mathbf{a}) = 0$.
	
	\begin{lemma}\label{lemma:volumeratio-GH}
		Let $(X,d, \H^N)$ be a $\mathrm{ncRCD}(K, N)$ space. There exists a function $\Psi_1(\varepsilon\mid K,N)$
		such that for every $x\in X$ and $r>0$ with $r^2K^-<\varepsilon$ where $K^-\coloneqq -K\vee 0$,
		\begin{itemize}
			\item  if $d_{\mathrm{GH}}(B_r^X(x),B_r^{\R^N}(0))\leq r\varepsilon$, then $\frac{\H^N(B_r^X(x))}{\omega_N
				r^N}\geq1 - \Psi_1(\varepsilon\mid K, N)$;
			\item conversely if $\frac{\H^N(B_r^X(x))}{\omega_N r^N}\geq 1-\varepsilon$, then
			$d_{\mathrm{GH}}(B_{r/2}^X(x), B_{r/2}^{\R^N}(0))\leq r\Psi_1(\varepsilon\mid K,N)$.
		\end{itemize}
	\end{lemma}
	
	\begin{proof}
		Denote by $X_r$ the rescaled space $(X,r^{-1}d,\H_{r^{-1}d}^N)$. Then $X_r$ is an $\mathrm{ncRCD}(r^2K,N)$
		space and $\H_{r^{-1}d}^N(B_1^{X_r}(x))=r^{-N}\H^N(B_r^X(x))$. Moreover,
		\[
		d_{\mathrm{GH}}(B_r^X(x), B_r^{\R^N}(0))=r\, d_{\mathrm{GH}}(B_1^{X_r}(x),B_1^{\R^N}(0)).
		\]
		Since $r^2K\geq-\varepsilon$, $X_r$ satisfies $\rcd(-\varepsilon,N)$. The first assertion therefore follows
		from the continuity of $\H^N$ under Gromov--Hausdorff convergence established in \cite[Theorem
		1.3]{DePhilippis-Gigli18}, applied through the standard compactness-contradiction argument. The second
		assertion follows in the same way from the volume rigidity theorem \cite[Theorem 1.6]{DePhilippis-Gigli18}.
		Scaling back to $(X,d,\H^N)$ concludes the proof.
	\end{proof}
	
	\begin{corollary}\label{cor:R-interior}
		There exists a function $\Psi(\varepsilon, \delta\mid K,N)$ s.t.\
		$\mathcal R_{\varepsilon,\delta}(X)\subset \operatorname{Int}(\mathcal R_{\Psi(\varepsilon,
			\delta),\delta/2}(X))$.
	\end{corollary}
	\begin{proof}
		Let $x\in\mathcal R_{\varepsilon,\delta}(X)$ and choose $t\in(\delta,2\delta)$ as in the definition.
		Applying the first part of \cref{lemma:volumeratio-GH} with $\varepsilon+4\delta^2K^-$, and then using
		continuity of ball volumes, we obtain an open neighborhood $U$ of $x$ such that
		\[
		\frac{\H^N(B_t^X(y))}{\omega_Nt^N}
		\geq1-2\Psi_1(\varepsilon+4\delta^2K^-\mid K,N),\qquad y\in U.
		\]
		Bishop--Gromov comparison with curvature lower bound $-K^-$ then gives
		\[
		\frac{\H^N(B_r^X(y))}{\omega_Nr^N}
		\geq\frac{\H^N(B_t^X(y))}{v_{-K^-,N}(t)}
		\frac{v_{-K^-,N}(r)}{\omega_Nr^N}\geq1-\eta,
		\qquad y\in U,\quad 0<r\leq t,
		\]
		where $\eta=\eta(\varepsilon,\delta\mid K,N)\to0$ as $(\varepsilon,\delta)\to(0,0)$.
		Here we use that $v_{-K^-,N}(r)/(\omega_Nr^N)\to1$ uniformly for $0<r\leq2\delta$ as $\delta\to0$.
		Enlarging $\eta$, we may assume that $\eta>4\delta^2K^-$.
		The second part of \cref{lemma:volumeratio-GH} then gives
		\[
		d_{\mathrm{GH}}(B_s^X(y),B_s^{\R^N}(0))\leq2s\Psi_1(\eta\mid K,N),
		\qquad y\in U,\quad 0<s\leq t/2.
		\]
		Set $\Psi(\varepsilon,\delta\mid K,N)\coloneqq2\Psi_1(\eta\mid K,N)$.
		Since $t/2>\delta/2$, we have $U\subset\mathcal R_{\Psi(\varepsilon,\delta\mid K,N),\delta/2}(X)$,
		which proves the claim.
	\end{proof}
	\begin{remark}\label{remark:R-interior}
		Let $f\colon X \to Y$ be a $1$-Lipschitz map between $\mathrm{ncRCD}(K, N)$ spaces such that
		$f_\#\H_X^N=\H_Y^N|_{f(X)}$. Then for all $\varepsilon,\delta>0$, $f(\mathcal R_{\varepsilon,\delta}(X))\subset \operatorname{Int}(\mathcal R_{\tilde \varepsilon,\delta/2}(Y))$ where
		$\tilde \varepsilon\coloneqq \Psi(\varepsilon,\delta\mid K,N)$ given as in \cref{cor:R-interior}.
		
		Indeed, for any $x\in \mathcal R_{\varepsilon, \delta}(X)$, by assumption we have
		\[
		\H_Y^N\bigl(B_r^Y(f(x))\bigr) \geq\H_Y^N\bigl(B_r^Y(f(x))\cap f(X)\bigr)
		=\H_X^N\bigl(f^{-1}(B_r^Y(f(x)))\bigr) \geq\H_X^N\bigl(B_r^X(x)\bigr)
		\]
		The conclusion follows by repeating the proof of \cref{cor:R-interior} on $Y$, starting from this volume
		comparison.
	\end{remark}
	\subsubsection*{Boundary of ncRCD spaces}
	We finally recall the notion of boundary for non-collapsed $\rcd$ spaces.
	Denote by $\operatorname{Tan}(X,x)$ the collection of tangent cones\footnote{Here a tangent cone means a
		pointed Gromov--Hausdorff limit of a blow-up sequence of $X$; see e.g. \cite[Definition
		2.6]{DePhilippis-Gigli18}. This notion should not be confused with the metric tangent introduced in
		\cref{sec:pre_rectifiable}.} of $X$ at $x$. For $0\leq k\leq N-1$, define the $k$-th singular stratum by
	\begin{equation}
		\mathcal S_X^k\coloneqq\left\{x\in X:\ \text{no tangent cone in }\operatorname{Tan}(X,x)\text{ splits off an
		}\R^{k+1}\text{-factor}\right\}.
	\end{equation}
	The singular strata satisfy
	\begin{equation}
		\mathcal S_X^0\subset\mathcal S_X^1\subset\cdots\subset\mathcal S_X^{N-1} =X\setminus\mathcal R(X),
	\end{equation}
	and
	\begin{equation}\label{ineq:dim_singular}
		\dim_{\mathrm{H}}\mathcal S_X^k\leq k,\qquad 0\leq k\leq N-1;
	\end{equation}
	see \cite[Theorem 1.8]{DePhilippis-Gigli18}.
	
	The reduced boundary is the top-dimensional singular stratum
	\begin{equation}
		\partial^*X\coloneqq\mathcal S_X^{N-1}\setminus\mathcal S_X^{N-2}.
	\end{equation}
	Following the boundary convention introduced in \cite{DePhilippis-Gigli18}, we define
	\begin{equation}
		\partial X\coloneqq\overline{\partial^*X} =\overline{\mathcal S_X^{N-1}\setminus\mathcal S_X^{N-2}},
	\end{equation}
	where the closure is taken in $X$. By \cite[Theorem~6.6]{Brue-Naber-Semola22}, $\partial X=\varnothing$ if
	and only if the inductively defined boundary $\mathcal F X$ of Kapovitch--Mondino \cite{Kapovitch-Mondino21} is empty. By
	\cite[Remark 6.4]{Brue-Naber-Semola22}, the Bishop--Gromov density satisfies $\theta_X(x)\leq 1/2$ on
	$\partial X$. In particular, $\partial X\cap\mathcal R(X)=\varnothing$.
	
	We define the \textbf{interior} of $X$ by $X^\circ\coloneqq X\setminus\partial X$. Moreover,
	\begin{equation}\label{eq:not_boundary_not_regular}
		X^\circ\setminus \mathcal R(X)=(X\setminus\mathcal R(X))\setminus\partial X\subset\mathcal S_X^{N-2}.
	\end{equation}
	Thus, due to \eqref{ineq:dim_singular}, the points which are neither regular points nor boundary points form
	a set of Hausdorff dimension at most $N-2$.
	\subsection{Local rigidity and length preservation}
	
	We first prove local isometry near regular points, using the results of
	\cref{sec:LVR_rectifiable1}. We then apply this local result to tangent maps
	to obtain a pointwise distance estimate and length preservation on all of $X$.
	Throughout this subsection and the next, $(X,d_X,\H_X^N)$ and $(Y,d_Y,\H_Y^N)$ are \(\mathrm{ncRCD}(K,N)\) spaces, and
	$f\colon X\to Y$ is $1$-Lipschitz and satisfies
	\[
	f_\#\H_X^N=\H_Y^N|_{f(X)}.
	\]
	By \cref{lemma:volume-preserving}, for every Borel set $A\subset X$ we have $\H_Y^N(f(A))=\H_X^N(A)$.
	Moreover, for every Borel set $B\subset Y$,
	\[
	\H_X^N(f^{-1}(B))=\H_Y^N(B\cap f(X))\leq\H_Y^N(B).
	\]
	These are the volume comparisons used below.
	
	With \(\Psi\) as in \cref{cor:R-interior}, chosen nondecreasing in each error variable, define
	\[
	\bar\varepsilon(\varepsilon) \coloneqq \max\bigl\{ \varepsilon,\Psi(\varepsilon,\varepsilon\mid K,N) \bigr\}.
	\]
	Then \(\bar\varepsilon\) is nondecreasing and
	\[
	\bar\varepsilon(\varepsilon)\geq\varepsilon,\qquad \bar\varepsilon(\varepsilon)\to0\quad\text{as
	}\varepsilon\searrow0,
	\]
	and
	\[
	\Psi(\varepsilon,\delta\mid K,N) \leq\bar\varepsilon(\varepsilon) \qquad\text{whenever
	}0<\delta\leq\varepsilon.
	\]
	
	\begin{lemma}[Local isometric embedding]
		\label{lem:regular-local-homeomorphism}
		Let \(\mathcal K_X\Subset X\) be compact and set
		\[
		m_{\mathcal K_X}\coloneqq\inf_{\mathcal K_X}\theta_X>0.
		\]
		Then there exists \(\varepsilon_{\mathcal K_X}>0\) such that, for every
		\(0<\varepsilon\leq\varepsilon_{\mathcal K_X}\), the following statements hold.
		\begin{enumerate}
			\item For every
			\(x\in\mathcal R_{\bar\varepsilon(\varepsilon)}(X)\), $\theta_X(x)>1-\frac{m_{\mathcal K_X}}2$. Moreover,
			\[
			f^{-1}(f(x))\cap\mathcal K_X=\{x\} \qquad\text{for every } x\in \mathcal R_{\bar\varepsilon(\varepsilon)}(X) \cap\mathcal K_X.
			\]
			
			\item If \(z\in\mathcal R_\varepsilon(X)\) and
			\(B_\tau(z)\Subset\mathcal K_X\) for some \(\tau>0\), then there exist open neighborhoods \(U_z\) of \(z\)
			and \(V_z\) of \(f(z)\) such that
			\[
			U_z\Subset \operatorname{Int}\bigl( \mathcal R_{\bar\varepsilon(\varepsilon)}(X) \bigr)\cap B_\tau(z),
			\qquad f(U_z)\subset V_z\Subset Y^\circ,
			\]
			\(f(U_z)\) is open in \(V_z\), and \(f|_{U_z}\colon U_z\to f(U_z)\) is an isometric embedding.
			
			\item If \(x\in\mathcal R(X)\), then $f(x)\in\mathcal R(Y)$ and $f^{-1}(f(x))=\{x\}$; if there is
			\(B_\tau(x)\Subset\mathcal K_X\) for some \(\tau>0\), then there exists \(c>0\) such that
			\[
			B_c(f(x))\subset f(\mathcal R_\varepsilon(X)).
			\]
		\end{enumerate}
	\end{lemma}
	
	\begin{proof}
		By \cref{lemma:volumeratio-GH}, choose \(\eta_{\mathcal K_X}>0\) so that, for every
		\(0<\eta\leq\eta_{\mathcal K_X}\), $\theta_X(x)>1-\frac{m_{\mathcal K_X}}2$ for every $x\in\mathcal R_\eta(X)$. Since \(\bar\varepsilon(\varepsilon)\to0\) as \(\varepsilon\searrow0\), we may choose
		\(\varepsilon_{\mathcal K_X}>0\) so that
		\begin{equation}\label{eq:epsilon_2}
			\bar\varepsilon(\varepsilon) \leq \min\left\{\eta_{\mathcal K_X},\frac{1}{1+K^-},\frac{\varepsilon_0(K,N)}2\right\},\quad \Psi_1(\bar\varepsilon(\varepsilon))<\frac12,
		\end{equation}
		for all $0<\varepsilon\leq\varepsilon_{\mathcal K_X}$, where $\Psi_1$ is taken from
		\cref{lemma:volumeratio-GH} and $\varepsilon_0(K,N)>0$ is the threshold in \cite[Theorem
		4.11(1)]{Kapovitch-Mondino21}.
		
		Suppose now that \(x,x'\in\mathcal K_X\), \(x\neq x'\), \(x\in\mathcal R_{\bar\varepsilon(\varepsilon)}(X)\), and \(f(x)=f(x')=y\). For \(0<r<d_X(x,x')/2\), the two balls are
		disjoint and, since \(f\) is \(1\)-Lipschitz,
		\[
		B_r(x)\sqcup B_r(x')\subset f^{-1}(B_r(y)).
		\]
		Measure preservation, followed by division by \(\omega_Nr^N\) and the limit \(r\searrow0\), gives
		\[
		1\geq\theta_Y(y) \geq\theta_X(x)+\theta_X(x') >1-\frac{m_{\mathcal K_X}}2+m_{\mathcal K_X}>1,
		\]
		a contradiction. This proves (1).
		
		For (2), let \(z\in\mathcal R_\varepsilon(X)\) and suppose that \(B_\tau(z)\Subset\mathcal K_X\). Put
		\(\varepsilon_1=\bar\varepsilon(\varepsilon)\). Since \(B_r(z)\subset f^{-1}(B_r(f(z)))\), measure
		preservation and (1) give
		\[
		\theta_Y(f(z))\geq\theta_X(z)>1-\frac{m_{\mathcal K_X}}2\geq\frac12.
		\]
		Thus \(f(z)\in Y^\circ\). By the definition of \(\mathcal R_\varepsilon(X)\), we may choose
		\(0<\delta_z\leq\varepsilon\) such that \(z\in\mathcal R_{\varepsilon,\delta_z}(X)\). Then
		\[
		\delta_z^2K^-<\varepsilon, \qquad \Psi(\varepsilon,\delta_z\mid K,N)\leq\varepsilon_1.
		\]
		By \cref{cor:R-interior},
		\[
		z\in \operatorname{Int}\bigl( \mathcal R_{\Psi(\varepsilon,\delta_z\mid K,N),\,\delta_z/2}(X) \bigr) \subset
		\operatorname{Int}\bigl(\mathcal R_{\varepsilon_1}(X)\bigr).
		\]
		Moreover, \cref{remark:R-interior} yields
		\[
		f(z)\in \operatorname{Int}\bigl( \mathcal R_{\Psi(\varepsilon,\delta_z\mid K,N),\,\delta_z/2}(Y) \bigr)
		\subset \operatorname{Int}\bigl(\mathcal R_{\varepsilon_1}(Y)\bigr).
		\]
		By \cite[Theorem 4.11(1)]{Kapovitch-Mondino21}, both \(\operatorname{Int}(\mathcal R_{\varepsilon_1}(X))\)
		and \(\operatorname{Int}(\mathcal R_{\varepsilon_1}(Y))\) are topological\linebreak[4]
		\(N\)-manifolds. We may therefore
		choose open \(N\)-manifold neighborhoods
		\[
		\begin{aligned}
			U_z &\Subset \operatorname{Int}(\mathcal R_{\varepsilon_1,\delta_z/2}(X))\cap B_\tau(z),\\
			V_z &\Subset \operatorname{Int}(\mathcal R_{\varepsilon_1,\delta_z/2}(Y))\cap Y^\circ
		\end{aligned}
		\]
		of \(z\) and \(f(z)\), respectively. After shrinking \(U_z\), continuity of \(f\) gives \(f(U_z)\subset
		V_z\). Statement (1) makes \(f|_{U_z}\) injective. The invariance of domain theorem therefore implies that
		\(f(U_z)\) is open in \(V_z\) and that \(f|_{U_z}\colon U_z\to f(U_z)\) is a homeomorphism.
		
		Now we show that $f\colon U_z\to f(U_z)$ is biLipschitz. To this aim, suppose otherwise. Then there exist
		distinct points $x_i,\tilde x_i\in U_z$ such that
		\begin{equation}
			d_Y(f(x_i), f(\tilde{x}_i)) \leq i^{-1} d_X(x_i, \tilde{x}_i).
		\end{equation}
		Let $d_i \coloneqq \frac{1}{2} d_X(x_i, \tilde{x}_i)$. It follows that $d_Y(f(x_i), f(\tilde{x}_i)) \leq
		\frac{2d_i}{i}$. Then, by \cref{lemma:volume-preserving},
		\begin{align}
			\H_Y^N\left(B_{(1+\frac{2}{i})d_i}(f(x_i))\right) &\geq \H_Y^N \left( B_{d_i}(f(\tilde{x}_i)) \cup
			B_{d_i}(f(x_i)) \right) \\
			&\geq  \H_Y^N\left( f\left( B_{d_i}(x_i)\sqcup B_{d_i}(\tilde{x}_i) \right) \right) \\
			&= \H_X^N( B_{d_i}(x_i)) +\H_X^N(B_{d_i}(\tilde{x}_i)).
		\end{align}
		Without loss of generality, we may assume $d_i\downarrow 0$ (otherwise by compactness of $\overline{U_z}$
		one can extract subsequences of $(x_i)$ and $(\tilde x_i)$ converging to $p$ and $q$ respectively s.t.
		$p\neq q$ but $f(p)=f(q)$, which contradicts (1)). Thus, we assume that $d_i<\delta_z/2$ for all $i$. By the
		Bishop--Gromov inequality and \cref{lemma:volumeratio-GH} (with the choice of $\varepsilon$ from
		\eqref{eq:epsilon_2}),
		\begin{align}
			d^{-N}_iv_{K,N}\big((1+2/i)d_i\big)&\geq d^{-N}_i\H_X^N( B_{d_i}(x_i))+d^{-N}_i\H_X^N( B_{d_i}(\tilde x_i))\\
			&\geq 2\omega_N(1-\Psi_1(\varepsilon_1))>\omega_N
		\end{align}
		for all $i$. But the left-hand side converges to $\omega_N$ as $i\to\infty$, giving a contradiction.
		
		By \cref{lemma:volume-preserving} and the injectivity of \(f|_{U_z}\), we have
		\((f|_{U_z})_\#(\H_X^N|_{U_z})=\H_Y^N|_{f(U_z)}\). Applying the local argument in the proof of
		\cref{thm:LVR-rectifiable} to \(f|_{U_z}\colon U_z\to Y\) shows that
		\(f|_{U_z}\) is a local isometry. Shrinking \(U_z\) around \(z\) proves (2).
		
		Finally, let \(x\) be as in (3), and put \(y=f(x)\). Since \(B_r(x)\subset f^{-1}(B_r(y))\), measure
		preservation gives \(\theta_Y(y)\geq\theta_X(x)=1\). Hence \(\theta_Y(y)=1\), so \(y\in\mathcal R(Y)\). If
		\(x'\neq x\) were another preimage of \(y\), the same disjoint-ball argument as in (1) would give
		\(1\geq\theta_Y(y)\geq1+\theta_X(x')>1\), a contradiction. Now suppose \(B_\tau(x)\Subset\mathcal K_X\).
		Choose \(0<s\leq\varepsilon_{\mathcal K_X}\) such that
		\[
		x\in\mathcal R_s(X), \qquad \bar\varepsilon(s)\leq\varepsilon.
		\]
		Applying (2) with the parameter \(s\), we obtain an open neighborhood \(U\subset\mathcal R_{\bar\varepsilon(s)}(X) \subset\mathcal R_\varepsilon(X)\) of \(x\) such that \(f(U)\) is open in \(Y\).
		Since \(y\in f(U)\), there exists \(c>0\) such that \(B_c(y)\subset f(U)\), proving (3).
	\end{proof}
	
	Recall that every tangent cone of an $\mathrm{ncRCD}(K,N)$ space is a metric cone,
	with the limiting basepoint as its vertex; see \cite[Proposition~2.8]{DePhilippis-Gigli18}.
	We write a metric cone over a base space $Z$ as $C(Z)$, and denote its metric by
	$d_C$ and its vertex by $o$. Each point other than $o$ has coordinates $(r,z)$,
	where $r=d_C(o,(r,z))>0$ and $z\in Z$. The radial projection is the map
	\[
	\pi\colon C(Z)\setminus\{o\}\to Z,\qquad \pi(r,z)=z.
	\]
	The following elementary fact will be used for maps between tangent cones.
	
	\begin{lemma}[Local geodesics from the vertex]\label{lem:cone-local-geodesic}
	    Let $(C(Z),d_C,o)$ be a metric cone with vertex $o$, and let $T>0$.
	    Suppose that $c\colon[0,T]\to C(Z)$ is continuous, $c(0)=o$, and
	    $c|_{(0,T)}$ is a unit-speed local geodesic. Then
	    \[
	    d_C(o,c(t))=t\qquad\text{for every }t\in[0,T].
	    \]
	\end{lemma}
	\begin{proof}
	    Set $q(t)\coloneqq d_C(o,c(t))^2$. Every point of $(0,T)$ lies in the interior
	    of a closed interval $I\subset(0,T)$ such that $c|_I$ is a unit-speed geodesic.
	    
	    If $c|_I$ avoids $o$, its radial projection has image $\pi(c(I))$, a geodesic
	    segment in $Z$ or a single point. The subcone $C(\pi(c(I)))$ is isometric to a Euclidean
	    sector, possibly reduced to a ray, in which $c|_I$ is a straight segment;
	    see \cite[I.5.10]{Bridson-Haefliger99}. Extend this segment to a line in the
	    Euclidean plane. Let $h$ be the distance
	    from $o$ to the line and $t_0$ the parameter at the perpendicular foot;
	    $t_0$ may lie outside $I$. Pythagoras gives
	    \begin{align}\label{eq:affine}
	    	q(t)-t^2=h^2+t_0^2-2t_0t \qquad\text{for }t\in I.
	    \end{align}
	    On the other hand, if $c(s)=o$ for some $s\in I$, then we have $q(t)=(t-s)^2$ for $t\in I$.
	    Thus, taking $h=0$ and $t_0=s$, we still have the identity \eqref{eq:affine}.
	    The coefficients agree on overlapping open intervals, so one can
	    extend this identity to $(0,T)$. Continuity at $0$ gives
	    $h^2+t_0^2=q(0)=0$. Therefore, $q(t)=t^2$ on $[0,T]$.
	\end{proof}
	
	We next use tangent cones to show that $f$ is asymptotically distance-preserving. Together with Lemma~\ref{lem:regular-local-homeomorphism}, the preceding lemma implies that the limiting map preserves distances from the basepoint. This asymptotic distance preservation then leads to length preservation, including for curves meeting the boundary.
	
	\begin{theorem}[Length preservation]\label{thm:length-isometry-RCD}
		Let $(X,d_X,\H_X^N)$ and $(Y,d_Y,\H_Y^N)$ be $\mathrm{ncRCD}(K,N)$ spaces,
		and let $f\colon X\to Y$ be $1$-Lipschitz with
		$f_\#\H_X^N=\H_Y^N|_{f(X)}$. Then, for every $x\in X$,
		\begin{equation}\label{eq:pointwise-distance-preservation}
		\lim_{\substack{z\to x\\z\neq x}}
		\frac{d_Y(f(z),f(x))}{d_X(z,x)}=1.
		\end{equation}
		In particular, \(\lip(f)(x)=1\) for every \(x\in X\) and $f$ is a length isometry.
	\end{theorem}
	\begin{proof}
	    Fix $x\in X$. Since $f$ is $1$-Lipschitz, failure of
	    \eqref{eq:pointwise-distance-preservation} would give $\varepsilon>0$ and
	    points $x_i\neq x$, $x_i\to x$, such that
	    \[
	    d_Y(f(x_i),f(x))\leq(1-\varepsilon)d_X(x_i,x).
	    \]
	    Set $r_i=d_X(x_i,x)$ and pass to a subsequence so that $r_i\downarrow0$.
	
	    Rescale the spaces and measures by setting
	    \[
	    X_i=(X,r_i^{-1}d_X,\mu_i,x),\qquad
	    Y_i=(Y,r_i^{-1}d_Y,\nu_i,f(x)),
	    \qquad \mu_i=r_i^{-N}\H_X^N,\quad \nu_i=r_i^{-N}\H_Y^N.
	    \]
	    The densities $\theta_X(x)$ and $\theta_Y(f(x))$ are positive and finite.
	    Pointed compactness and \cite[Theorem~1.2 and Proposition~2.8]{DePhilippis-Gigli18}
	    therefore give, after passing to a subsequence, pointed measured
	    Gromov--Hausdorff limits
	    \[
	    X_i\longrightarrow(C_X,d_{C_X},\mu,o_X),\qquad
	    Y_i\longrightarrow(C_Y,d_{C_Y},\nu,o_Y),
	    \]
	    where $C_X=C(Z_X)$ and $C_Y=C(Z_Y)$ are $\mathrm{ncRCD}(0,N)$ metric cones,
	    equipped with $\mu=\H_{C_X}^N$, and $\nu=\H_{C_Y}^N$, respectively.
	    The rescaled maps $f_i\colon X_i\to Y_i$ are $1$-Lipschitz and preserve
	    the basepoints. A diagonal Arzel\`a--Ascoli argument gives a further
	    subsequence converging locally uniformly to a $1$-Lipschitz map
	    $F\colon C_X\to C_Y$ with $F(o_X)=o_Y$.
	    In particular, $f_i(z_i)\to F(z)$ whenever $z_i\to z$ under realizations
	    of the pointed convergences.
	
		By local convergence of the maps and measures and exhaustion of $C_X$,
		the inequalities $(f_i)_\#\mu_i\leq\nu_i$ pass to the limit and give
		$F_\#\mu\leq\nu$.
		
		The cone $C_X$ is proper, so its image $S\coloneqq F(C_X)$ is a countable union
		of compact sets and hence Borel. For every Borel set $B\subset C_Y$, the
		$1$-Lipschitz property of $F$ and the inequality just proved give
		\[
		\nu(B\cap S)=\H_{C_Y}^N\bigl(F(F^{-1}(B))\bigr)
		\leq\mu(F^{-1}(B))=(F_\#\mu)(B\cap S)\leq\nu(B\cap S).
		\]
		Since $F_\#\mu$ is concentrated on $S$, it follows that
		$F_\#\mu=\nu|_S$.
	
	    Now fix $p\in\mathcal R(C_X)\setminus\{o_X\}$ and let $\rho\colon[0,L]\to C_X$
	    be the unit-speed radial segment from $o_X$ to $p$.
	    Cone dilations preserve regular points, so $\rho(t)\in\mathcal R(C_X)$
	    for every $t\in(0,L]$. By the measure identity above and
	    Lemma~\ref{lem:regular-local-homeomorphism}\,(2), $F$ is a local isometry at these points.
	    Thus $F\circ\rho$ is a unit-speed local geodesic on $(0,L)$.
	    Since $F(\rho(0))=o_Y$, Lemma~\ref{lem:cone-local-geodesic} gives
	    $d_{C_Y}(o_Y,F(p))=L$. By density of the regular set and continuity of $F$,
	    \[
	    d_{C_Y}(o_Y,F(p))=d_{C_X}(o_X,p)\qquad\text{for every }p\in C_X.
	    \]
	
	    Since $d_{X_i}(x,x_i)=1$, the pointed convergence gives, after a further subsequence,
	    $x_i\to p\in C_X$ with $d_{C_X}(o_X,p)=1$. Hence
	    \[
	    1=d_{C_Y}(o_Y,F(p))
	    =\lim_{i\to\infty}r_i^{-1}d_Y(f(x),f(x_i))\leq1-\varepsilon,
	    \]
	    a contradiction. This proves \eqref{eq:pointwise-distance-preservation}.
	
		Finally, let $\gamma$ be a rectifiable curve of positive length and parametrize it by arc length
		on $[0,L]$; the zero-length case is immediate. Both $\gamma$ and
		$\alpha\coloneqq f\circ\gamma$ are Lipschitz. Applying
		\eqref{eq:pointwise-distance-preservation} at $\gamma(t)$ gives
		$|\dot\alpha_t|=|\dot\gamma_t|$ for almost every $t$. The length formula yields
		\[
		L_Y(f\circ\gamma)=\int_0^L|\dot\alpha_t|\,\d t
		=\int_0^L|\dot\gamma_t|\,\d t=L_X(\gamma),
		\]
		which completes the proof.
	\end{proof}
	
	\subsection{Openness and injectivity on the interior}
	
	We now use the pointwise distance estimate to prove that $f$ is open on
	$X^\circ$. The other ingredient is the following dimension-comparison lemma,
	which prevents a set of codimension at least two from separating a point
	from a set of positive volume along a bounded family of geodesics.
	\begin{lemma}[Dimension comparison for a bounded family of geodesics]
		\label{lem:localized-dimension-comparison}
		Let $(Z,d_Z,\H_Z^N)$ be a $\mathrm{ncRCD}(K,N)$ space,
		with no assumption on $\partial Z$, and let
		$\mathcal K_Z\subset Z$ be compact.
		Fix $z_0\in\mathcal K_Z$ and an $\H_Z^N$-measurable set
		$\Omega\subset\mathcal K_Z$ with $\H_Z^N(\Omega)>0$.
		For each $z\in\Omega$, choose a minimizing geodesic
		$\gamma_z$ from $z_0$ to $z$, contained in $\mathcal K_Z$,
		and a point $w_z\in\gamma_z$.
		If
		\[
		\Sigma:=\{w_z:z\in\Omega\},
		\qquad \operatorname{dist}_Z(z_0,\Sigma)>0,
		\]
		then
		\[
		\dim_{\mathrm H}\Sigma\geq N-1.
		\]
	\end{lemma}
	
	\begin{proof}
		We use a localized version of the covering argument in \cite[Lemma 5.2 and Corollary 5.3]{CDNPSG24}. Write
		\(\m=\H_Z^N\), \(D=\operatorname{diam}(\mathcal K_Z)\), \(a=\operatorname{dist}_Z(z_0,\Sigma)>0\), and
		\(b=a/(2D)\). Parameterize the chosen geodesics on \([0,1]\) with constant speed. By the almost-everywhere
		uniqueness of geodesics from a fixed point \cite[Corollary 1.4]{Gigli-Rajala-Sturm16}, after discarding a
		null set of endpoints, the maps \(F_t(z)\coloneqq \gamma_z(t)\) are measurable for all $t\in [b,1]$. Since
		each $\gamma_z$ is continuous, the map $(t,z)\mapsto F_t(z)$ is jointly measurable. By the
		MCP condition \cite[Corollary 1.5]{Gigli-Rajala-Sturm16}, there exists \(C\) depending only
		on \(K,N,D,b\) s.t.
		\begin{equation}\label{ineq:densitybound}
			(F_t)_\#(\m|_\Omega)\leq C\,\m \qquad\text{for every }t\in[b,1].
		\end{equation}
		In particular, these facts apply to the prescribed geodesics, and require no boundary assumption.
		
		Let \(\{B_{r_j}(q_j)\}_j\) be a countable cover of $\Sigma$ with \(0<r_j<a/2\). For each remaining
		\(z\in\Omega\), choose \(j\) such that \(w_z=\gamma_z(t_z)\in B_{r_j}(q_j)\). Since \(t_z\geq a/D=2b\) and
		\(\gamma_z\) has speed at most \(D\), we have
		\[
		\gamma_z([t_z-r_j/D,t_z])\subset B_{2r_j}(q_j), \qquad J_z\coloneqq [t_z-r_j/D,t_z]\subset[b,1].
		\]
		Hence for all those $z\in \Omega$,
		\begin{equation}
			\frac{r_j}{D}=\L^1(J_z)\leq \int^1_b \mathds{1}_{B_{2r_j}(q_j)}(F_t(z))\d t\,.
		\end{equation}
		Consequently, Fubini's theorem and the density bound \eqref{ineq:densitybound} give
		\[
		0<\m(\Omega) \leq\sum_j\frac{D}{r_j} \int_b^1\int_\Omega \mathds{1}_{B_{2r_j}(q_j)}(F_t(z))\,\d\m(z)\,\d t
		\leq CD\sum_j r_j^{-1}\m(B_{2r_j}(q_j)).
		\]
		By Bishop--Gromov comparison and \(\theta_Z\leq1\), for sufficiently small \(r_j\) the last expression is
		bounded by \(C'\sum_j r_j^{N-1}\), with \(C'\) independent of the cover. Taking the infimum over such covers
		yields \(\H^{N-1}_Z(\Sigma)>0\), and hence \(\dim_{\mathrm{H}}\Sigma\geq N-1\).
	\end{proof}
	
	\begin{lemma}[Openness]\label{lemma:openness}
	    The map $f|_{X^\circ}\colon X^\circ\to Y$ is open. In particular, $f(X^\circ)$ is open in $Y$.
	\end{lemma}
	\begin{proof}
	    Let $V\subset X^\circ$ be open, fix $x\in V$, and write $a=f(x)$.
	    By \eqref{eq:pointwise-distance-preservation}, we can choose $r>0$ so small that
	    $\mathcal K_X\coloneqq\overline{B_r(x)}\subset V$ and
	    \[
	    d_Y(f(z),a)\geq\tfrac12 d_X(z,x)\qquad\text{for every }z\in\mathcal K_X.
	    \]
	    Set $A\coloneqq f(\mathcal K_X)$, which is compact. In this proof,
	    $\partial_Y A$ denotes the topological boundary of $A$ in $Y$.
	    We claim that
	    \begin{equation}\label{eq:compact-image-boundary}
	        \partial_Y A\cap B_{r/2}(a)\subset f(\mathcal K_X\setminus\mathcal R(X)).
	    \end{equation}
	    Indeed, if $y$ belongs to the set on the left, compactness gives a point
	    $z\in\mathcal K_X$ with $f(z)=y$, and the preceding estimate gives
	    $d_X(z,x)\leq2d_Y(y,a)<r$. If $z$ were regular,
	    Lemma~\ref{lem:regular-local-homeomorphism}\,(2), applied in a compact neighborhood
	    of $z$, would give an open neighborhood $W\subset B_r(x)$ of $z$ such that
	    $f(W)$ is open in $Y$. Then $y\in f(W)\subset A$, contradicting
	    $y\in\partial_Y A$. This proves \eqref{eq:compact-image-boundary}.
	    Since $\mathcal K_X\subset X^\circ$, \eqref{eq:not_boundary_not_regular},
	    \eqref{ineq:dim_singular}, and the Lipschitz property of $f$ imply
	    \[
	    \dim_{\mathrm H}\bigl(\partial_Y A\cap B_{r/2}(a)\bigr)\leq N-2.
	    \]
	
	    Choose $z_0\in B_{r/8}(x)\cap\mathcal R(X)$ and set $y_0=f(z_0)$.
	    The same application of Lemma~\ref{lem:regular-local-homeomorphism}\,(2) shows that
	    $y_0$ lies in the interior of $A$. Also, $d_Y(y_0,a)<r/8$.
	    Suppose that $\Omega\coloneqq B_{r/8}(a)\setminus A$ is nonempty.
	    Since $A$ is closed, $\Omega$ is open and has positive $\H_Y^N$-measure.
	    For each $y\in\Omega$, choose a minimizing geodesic from $y_0$ to $y$.
	    Its length is less than $r/4$, so it is contained in $B_{3r/8}(a)$.
	    It must meet $\partial_Y A$ at some point $w_y$. Thus
	    \[
	    \Sigma\coloneqq\{w_y:y\in\Omega\}
	    \subset\partial_Y A\cap B_{r/2}(a),\qquad
	    \operatorname{dist}_Y(y_0,\Sigma)>0,
	    \]
	    where the positive distance follows because $y_0$ is an interior point of $A$.
	    Lemma~\ref{lem:localized-dimension-comparison}, applied in
	    $\mathcal K_Y=\overline{B_{r/2}(a)}$, now gives
	    \[
	    N-1\leq\dim_{\mathrm H}\Sigma\leq N-2,
	    \]
	    a contradiction. Hence $B_{r/8}(a)\subset A\subset f(V)$.
	    Since $x\in V$ was arbitrary, $f(V)$ is open.
	\end{proof}
	
	Openness and volume preservation now give injectivity on the interior.
	The argument also rules out a boundary point with the same image as an interior point.
	
	\begin{lemma}[Interior injectivity]\label{lemma:int_injective}
	    For each $x\in X^\circ$, $f^{-1}(f(x))=\{x\}$. In particular, the restriction
	    $f|_{X^\circ}\colon X^\circ\to f(X^\circ)$ is injective.
	\end{lemma}
	\begin{proof}
	    Suppose that $p\in X^\circ$ and $q\in X\setminus\{p\}$ satisfy $f(p)=f(q)$.
	    Choose a bounded open neighborhood $V$ of $p$ with
	    $\overline V\subset X^\circ$ and $q\notin\overline V$.
	    By Lemma~\ref{lemma:openness}, $f(V)$ is an open neighborhood of $f(q)$.
	    Continuity of $f$ therefore gives a ball $W$ about $q$ such that
	    $W\cap V=\varnothing$ and $f(W)\subset f(V)$.
	    By Lemma~\ref{lemma:volume-preserving},
	    \[
	    \H_X^N(V)+\H_X^N(W)
	    =\H_Y^N(f(V\cup W))
	    =\H_Y^N(f(V))
	    =\H_X^N(V).
	    \]
	    This contradicts $\H_X^N(V)<\infty$ and $\H_X^N(W)>0$.
	\end{proof}
	
	\subsection{Proof of the main theorems}
	Let $(Z,d_Z)$ be a metric space. For $A\subset Z$ a Borel subset, we define the intrinsic (extended) metric
	$d_{A}^\ell$ on $A$ by
	\[
	d_{A}^\ell(x,y)\coloneqq \inf \{L_Z(\gamma):\gamma\in C([0,1];A),\gamma(0)=x,\gamma(1)=y\}
	\]
	where we set $\inf \varnothing=+\infty$.
	
	\begin{lemma}\label{lem:interior-ambient-metric}
		Let \(Z=X,Y\) be non-collapsed $\rcd(K,N)$ spaces. The interior \(Z^\circ\) is path-connected and
		\[
		d_{Z^\circ}^\ell(z_0,z_1)=d_Z(z_0,z_1) \qquad\text{for all }z_0,z_1\in Z^\circ.
		\]
	\end{lemma}
	\begin{proof}
		Since \(Z\) is non-collapsed, its essential dimension is \(N\), and its \(N\)-regular set is \(\mathcal R(Z)\). By \cite[Theorem 6.5]{Deng25}, \(\mathcal R(Z)\) is weakly convex. Thus, for every \(u,v\in\mathcal R(Z)\) and every \(\eta>0\), there is a piecewise geodesic curve \(\alpha\subset\mathcal R(Z)\) from \(u\)
		to \(v\) such that
		\[
		L_Z(\alpha)\leq d_Z(u,v)+\eta.
		\]
		As observed above, \(\mathcal R(Z)\subset Z^\circ\).
		
		Fix \(z_0,z_1\in Z^\circ\). The regular set has full \(\H_Z^N\)-measure and is therefore dense in \(Z\).
		Choose
		\[
		u_j,v_j\in\mathcal R(Z),\qquad u_j\to z_0,\quad v_j\to z_1,
		\]
		so close to the respective endpoints that minimizing geodesics \([z_0u_j]\) and \([v_jz_1]\) are contained
		in \(Z^\circ\). Applying weak convexity with \(\eta=1/j\) and concatenating the resulting three curves gives
		\[
		d_{Z^\circ}^\ell(z_0,z_1) \leq d_Z(z_0,u_j)+d_Z(u_j,v_j)+\frac1j+d_Z(v_j,z_1).
		\]
		Letting \(j\to\infty\) yields \(d_{Z^\circ}^\ell(z_0,z_1)\leq d_Z(z_0,z_1)\). The reverse inequality follows
		from the definition of length. In particular, \(d_{Z^\circ}^\ell\) is finite on \(Z^\circ\), so \(Z^\circ\)
		is path-connected.
	\end{proof}

	\begin{remark}
		Let \((X,d_X,\H_X^N)\), \((Y,d_Y,\H_Y^N)\) be \(\mathrm{ncRCD}(K,N)\) spaces. Assume also that
		\[
		\inf_{x\in X}\H_X^N(B_1(x))>0.
		\]
		Let \(f\colon X\to Y\) be a \(1\)-Lipschitz map such that \(f_\#\H_X^N=\H_Y^N\). Then \(f\) is surjective.
		
		To prove it, fix an arbitrary point \(y\in Y\). The identity \(f_\#\H_X^N=\H_Y^N\), together with the full
		support of \(\H_Y^N\), implies that \(f(X)\) is dense in \(Y\), thus we can find a sequence
		\((x_n)_n\subset X\) such that \(d_Y(f(x_n),y)<\frac{1}{n}\). If \((x_n)_n\) is an unbounded sequence, then
		we can extract a subsequence \((x_{n_i})_i\) such that \(d_X(x_{n_i},x_{n_j})\geq 2\) for every
		\(i,j\in\N\) with \(i\neq j\), so that the balls \(\{B_1(x_{n_i}):i\in\N\}\) are pairwise disjoint and thus
		accordingly, letting \(v\coloneqq\inf_{x\in X}\H_X^N(B_1(x))>0\), we have that
		\[
		\begin{split}
			+\infty&=\sum_{i=1}^\infty v\leq\sum_{i=1}^\infty\H_X^N(B_1(x_{n_i}))=\H_X^N\bigg(\bigcup_{i=1}^\infty
			B_1(x_{n_i})\bigg) \leq\H_X^N\bigg(\bigcup_{i=1}^\infty f^{-1}(f(B_1(x_{n_i})))\bigg)\\
			&\leq\H_X^N\bigg(\bigcup_{i=1}^\infty
			f^{-1}(B_1(f(x_{n_i})))\bigg)\leq\H_X^N(f^{-1}(B_2(y)))=\H_Y^N(B_2(y))<+\infty,
		\end{split}
		\]
		which leads to a contradiction. Therefore, \((x_n)_n\) is a bounded sequence. Since \(X\) is proper, up to
		a subsequence we have that \(x_n\) converges to some point \(x\in X\), thus accordingly \(y=\lim_n
		f(x_n)=f(x)\in f(X)\). This shows that \(f\) is surjective. The full pushforward identity is essential
		here: the inclusion $\R^{N-1}\times[0,\infty)\hookrightarrow\R^N$ satisfies $f_\#\H_X^N=\H_Y^N|_{f(X)}$ and
		the uniform lower volume bound, but is not surjective.
	\end{remark}

	\begin{theorem}\label{thm:LVR_RCD}
		Let \((X,d_X,\H_X^N)\), \((Y,d_Y,\H_Y^N)\) be \(\mathrm{ncRCD}(K,N)\) spaces. Let \(f\colon X\to Y\) be a
		\(1\)-Lipschitz map such that \(f_\#\H_X^N=\H_Y^N|_{f(X)}\). Then $f(X^\circ)\subset Y^\circ$, and the restriction
		\[
		f^\circ\coloneqq f|_{X^\circ}\colon (X^\circ, d_X) \to (f(X^\circ), d_{f(X^\circ)}^\ell)
		\]
		is an isometry. Moreover, if in addition we assume that $\mathcal R(Y)\subset f(X^\circ)$, then $f$ is
		surjective and $f\colon (X,d_X)\to (Y,d_Y)$ is an isometry.
	\end{theorem}
	The condition $\mathcal R(Y)\subset f(X^\circ)$ holds if $f$ is surjective and $f(\partial X)\subset\partial
	Y$. This boundary condition is automatic for bi-Lipschitz homeomorphisms. Indeed, \cite[Theorem
	1.6]{Kapovitch-Mondino21} gives $f(\mathcal F X)=\mathcal F Y$ for the inductively defined boundaries. Since
	$\overline{\mathcal F X}=\partial X$ and $\overline{\mathcal F Y}=\partial Y$ by \cite[Theorem
	6.6(i)]{Brue-Naber-Semola22}, taking closures gives $f(\partial X)=\partial Y$.
	
	\begin{proof}
	    Set $U\coloneqq f(X^\circ)$. Since $f$ is $1$-Lipschitz,
	    Lemma~\ref{lem:interior-ambient-metric} gives
	    \begin{equation}\label{eq:1Lipschitz}
	        d_U^\ell(f(x),f(x'))\leq d_{X^\circ}^\ell(x,x')=d_X(x,x')
	        \qquad\text{for all }x,x'\in X^\circ.
	    \end{equation}
	    By Lemmas~\ref{lemma:openness} and~\ref{lemma:int_injective},
	    $f^\circ$ is a homeomorphism onto the open set $U$.
	    Write $g\colon U\to X^\circ$ for its continuous inverse.
	
	    Let $\sigma\colon[0,L]\to U$ be a rectifiable curve parametrized by arc length,
	    and set $\gamma\coloneqq g\circ\sigma$.
	    We first show that $\gamma$ is $1$-Lipschitz. For each $t\in[0,L]$ and each
	    $\varepsilon>0$, continuity of $\gamma$ and
	    \eqref{eq:pointwise-distance-preservation} give
	    \begin{equation}\label{eq:inverse-curve-local-bound}
	        d_X(\gamma(s),\gamma(t))
	        \leq(1+\varepsilon)d_Y(\sigma(s),\sigma(t))
	        \leq(1+\varepsilon)|s-t|
	    \end{equation}
	    for all $s$ sufficiently close to $t$.
	    To pass from this local bound to all pairs of parameters, fix $0\leq a<b\leq L$.
	    Let $T$ be the supremum of the times $t\in[a,b]$ for which
	    \[
	    d_X(\gamma(a),\gamma(u))\leq(1+\varepsilon)(u-a)
	    \qquad\text{for every }u\in[a,t].
	    \]
	    The local bound at $a$ gives $T>a$, and continuity gives the same inequality
	    for every $u\in[a,T]$. If $T<b$, the local bound at $T$ and the triangle
	    inequality extend it to an interval beyond $T$, a contradiction.
	    Thus $T=b$. Letting $\varepsilon\downarrow0$ proves that $\gamma$ is $1$-Lipschitz.
	    Consequently,
	    \[
	    d_X(g(\sigma(0)),g(\sigma(L)))\leq L=L_Y(\sigma).
	    \]
	    Taking the infimum over curves joining any two points of $U$ shows that
	    $g\colon(U,d_U^\ell)\to(X^\circ,d_X)$ is $1$-Lipschitz.
	    Together with \eqref{eq:1Lipschitz}, this proves that $f^\circ$ is an isometry.
        Since the intrinsic metric of an open subset agrees locally with the ambient metric, $f|_{X^\circ}$ is locally an isometry between open subsets of $X$ and $Y$, which preserves the boundary. Hence $f(X^\circ)\subseteq Y^\circ$.
        
	    Assume now that $U\supset\mathcal R(Y)$. By \cite[Theorem~6.5]{Deng25},
	    $d_U^\ell=d_Y$ on $\mathcal R(Y)$.
	    Since $f(\mathcal R(X))\subset\mathcal R(Y)$ by
	    Lemma~\ref{lem:regular-local-homeomorphism}\,(3), we obtain
	    \[
	    d_X(x,x')=d_U^\ell(f(x),f(x'))=d_Y(f(x),f(x'))
	    \qquad\text{for all }x,x'\in\mathcal R(X).
	    \]
	    Density of $\mathcal R(X)$ and continuity of $f$ extend this identity to all
	    $x,x'\in X$. Hence $f$ is an isometric embedding. Its image is closed because
	    $X$ is complete, and it contains the dense set $\mathcal R(Y)$. Thus $f(X)=Y$.
	\end{proof}
	
	\begin{remark}[Two-dimensional case]\label{remark:lowerdim}
	    For $N=2$, \cite{Lytchak-Stadler23} proves that every
	    $\mathrm{ncRCD}(K,2)$ space is an Alexandrov space with curvature bounded below by $K$.
	    In the finite-volume setting, the interior rigidity in Theorem~\ref{thm:LVR_RCD}
	    and the length preservation in Theorem~\ref{thm:length-isometry-RCD}
	    are already contained in \cite[Theorem~A]{Li15}.
	\end{remark}
	\begin{theorem}\label{thm:LVR_RCD_2}
		Let \((X,d_X,\H_X^N)\), \((Y,d_Y,\H_Y^N)\) be \(\mathrm{ncRCD}(K,N)\) spaces. Let \(f\colon X\to Y\) be a
		\(1\)-Lipschitz map such that \(f_\#\H_X^N=\H_Y^N|_{f(X)}\). Then the following conditions are equivalent:
		\begin{itemize}
			\item[\(\rm i)\)] \(f(\partial X)\subseteq\partial Y\).
			\item[\(\rm ii)\)] \(f(X^\circ)=Y^\circ\) and \(f(\partial X)=\partial Y\).
			\item[\(\rm iii)\)] \(f\colon X\to Y\) is an isometric bijection.
		\end{itemize}
	\end{theorem}
    
	\begin{proof}
		The implication \({\rm iii)}\Rightarrow{\rm ii)}\) can be easily checked, while the implication \({\rm
			ii)}\Rightarrow{\rm i)}\) is trivial. It remains to prove that \({\rm i)}\Rightarrow{\rm iii)}\). To this
		aim, our plan is to show that \({\rm i)}\) implies that \(Y^\circ\subset f(X^\circ)\), thus in particular
		\(\mathcal R(Y)\subset f(X^\circ)\), whence iii) follows thanks to the last assertion of Theorem
		\ref{thm:LVR_RCD}. We argue by contradiction: suppose there exists \(y\in Y^\circ\setminus f(X^\circ)\).
		Choose $x_0\in\mathcal R(X)$. By \cref{lem:regular-local-homeomorphism}~(3), $f(x_0)\in\mathcal R(Y)\subset
		Y^\circ$. Lemma \ref{lem:interior-ambient-metric} ensures that \(Y^\circ\) is connected by rectifiable
		paths, thus there exists a constant-speed rectifiable curve \(\sigma\colon[0,1]\to Y^\circ\) such that
		\(\sigma_0=f(x_0)\) and \(\sigma_1=y\). Define \(\delta\coloneqq\inf\{t\in[0,1]:\sigma_t\notin
		f(X^\circ)\}\). Since \(f(X^\circ)\) is open by Lemma \ref{lemma:openness}, we have \(\delta>0\) and
		\(\sigma_\delta\notin f(X^\circ)\). The map \(f|_{X^\circ}\colon X^\circ\to f(X^\circ)\) is bijective by
		Theorem \ref{thm:LVR_RCD}, thus it makes sense to define \(\gamma\colon[0,\delta)\to X^\circ\) as
		\(\gamma_t\coloneqq f^{-1}(\sigma_t)\) for every \(t\in[0,\delta)\). For any \(s,t\in[0,\delta)\) with
		\(s<t\), we have that
		\[
		d_X(\gamma_s,\gamma_t)=d_{f(X^\circ)}^\ell(\sigma_s,\sigma_t)\leq L_Y(\sigma|_{[s,t]})=|t-s|L_Y(\sigma),
		\]
		where the first identity is guaranteed by the fact that \(f|_{X^\circ}\colon(X^\circ, d_X)\to(f(X^\circ),
		d_{f(X^\circ)}^\ell)\) is an isometry, as we proved in Theorem \ref{thm:LVR_RCD}. In particular, the limit
		\(x_1\coloneqq\lim_{t\nearrow\delta}\gamma_t\in X\) exists due to the completeness of \(X\). Then
		\(f(x_1)=\lim_{t\nearrow\delta}f(\gamma_t)=\lim_{t\nearrow\delta}\sigma_t=\sigma_\delta\notin f(X^\circ)\),
		so that necessarily \(x_1\in\partial X\) and thus \(\sigma_\delta=f(x_1)\in f(\partial X)\subset\partial
		Y\), which leads to a contradiction with the fact that \(\sigma\) is a curve in \(Y^\circ\). Therefore, we
		have \(Y^\circ\setminus f(X^\circ)=\varnothing\).
	\end{proof}

    Theorem~\ref{thm:length-isometry-RCD}, Lemma~\ref{lemma:openness}, and
	Theorems~\ref{thm:LVR_RCD} and~\ref{thm:LVR_RCD_2} together complete the proof of
	Theorem~\ref{thm:main_2}.

	\subsection{Almost rigidity}
    Using the convergence of boundary measures established in \cite[Theorem~1.8]{Brue-Naber-Semola22}, we obtain an almost-rigidity result (Theorem \ref{thm:main_almost}).
    It extends \cite[Corollary~0.1]{Li-Wang14} from non-collapsed Ricci limit spaces to ncRCD spaces, allowing non-empty boundary under an additional condition on the boundary volumes.

	Let $(X,d,\H^N)$ be a \(\mathrm{ncRCD}(K,N)\) space with non-empty boundary.
	Recall that, at every $p\in \mathcal S_X^{N-1}\setminus\mathcal S_X^{N-2}$, the tangent cone at $p$ is unique and isomorphic to $\R^{N}_+$; see  \cite[Theorem~1.4~(iii)]{Brue-Naber-Semola22}.
	By the $\varepsilon$-regularity \cite[Theorem~1.2]{Brue-Naber-Semola22}, $\H^{N-1}(B_r(p)\cap \partial X)>0$ for every $r>0$.
	Together with \cite[Theorem~1.4~(i)]{Brue-Naber-Semola22}, this shows that $\partial X$ is $(N-1)$-rectifiable in the sense of \cref{sec:pre_rectifiable}\footnote{Compared with the convention in
	\cite{Brue-Naber-Semola22}, our definition additionally requires
	positive Hausdorff measure.}.
    Moreover, $\H^{N-1}|_{\partial X}$ is a nonzero, locally finite measure satisfying
	\begin{equation}\label{eq:support-boundary-measure}
		\supp (\H^{N-1}|_{\partial X})=\partial X.
	\end{equation}

Before proving the almost rigidity theorem, we give a further criterion for a surjective volume-preserving $1$-Lipschitz map to be an isometry, based on the comparison of boundary measures.
\begin{proposition}\label{prop:LVR_RCD_3}
	Let $(X,d_X,\H_X^N)$ and $(Y,d_Y,\H_Y^N)$ be
	${\rm ncRCD}(K,N)$ spaces, and let $f\colon X\to Y$ be a surjective
	$1$-Lipschitz map satisfying $f_\#\H_X^N=\H_Y^N$.
	Then the following conditions are equivalent:
	\begin{itemize}
		\item[\textup{(i)}] $f_\#(\H_X^{N-1}|_{\partial X})\leq \H_Y^{N-1}|_{\partial Y}$.
		\item[\textup{(ii)}] $f\colon X\to Y$ is an isometric bijection.
	\end{itemize}
If in addition $\H_Y^{N-1}(\partial Y)<\infty$, then these conditions are also equivalent to
	\begin{itemize}
	\item[\textup{(iii)}] $\H_X^{N-1}(\partial X)\leq\H_Y^{N-1}(\partial Y)$.
	\end{itemize}
\end{proposition}
\begin{proof}
By \cref{thm:LVR_RCD}, $f(X^\circ)\subseteq Y^\circ$.
Surjectivity therefore gives $\partial Y\subseteq f(\partial X)$.
Since $f$ is $1$-Lipschitz, for every Borel set $B\subseteq Y$,
\begin{equation}
	\H_Y^{N-1}(B\cap\partial Y)\leq \H_Y^{N-1}\big(B\cap f(\partial X)\big)\leq\H_X^{N-1}\big(f^{-1}(B)\cap\partial X\big).
\end{equation}
In other words, $\H_Y^{N-1}|_{\partial Y}\leq f_\#(\H_X^{N-1}|_{\partial X})$.

Assume (i). Then $\H_Y^{N-1}|_{\partial Y}= f_\#(\H_X^{N-1}|_{\partial X})$.
Thus
\begin{equation}
 \H_X^{N-1}|_{\partial X} \big( f^{-1}(Y^\circ) \big)=f_\#(\H_X^{N-1}|_{\partial X})(Y^\circ)=\H_Y^{N-1}|_{\partial Y}(Y^\circ)=0.
\end{equation}
Since $f^{-1}(Y^\circ)$ is open, \eqref{eq:support-boundary-measure} implies that $\partial X\cap f^{-1}(Y^\circ)=\emptyset$.
Therefore $f(\partial X)\subseteq\partial Y$, and by \cref{thm:LVR_RCD_2}, $f$ is an isometry.
The implication (ii)$\Rightarrow$(i) is immediate.

Assume now (iii) with $\H_Y^{N-1}(\partial Y)<\infty$.
Then the two finite measures $\H_Y^{N-1}|_{\partial Y}$ and $f_\#(\H_X^{N-1}|_{\partial X})$ have equal total mass, so the comparison $\H_Y^{N-1}|_{\partial Y}\leq f_\#(\H_X^{N-1}|_{\partial X})$ forces them to be equal.
In particular, (i) holds, completing the proof.
\end{proof}
    
	\begin{theorem}[Almost rigidity]\label{thm:almost_rigidity}
		Fix $K\in\R$, $N\in\N$, and $D,v>0$.
		For every $\varepsilon>0$, there exists $\delta=\delta(K,N,D,v,\varepsilon)>0$ such that the following holds.
		
		Let $(X,d_X,\H_X^N)$ and $(Y,d_Y,\H_Y^N)$ be ${\rm ncRCD}(K,N)$ spaces satisfying
		\begin{equation}\label{eq:assumptions_almost_LVR}
			\diam(X)\leq D,
			\qquad
			\H_Y^N(Y)\geq v.
		\end{equation}
		Suppose that a map $f\colon X\to Y$ satisfies
		\begin{equation}\label{eq:defect_almost_LVR}
			\begin{gathered}
				d_Y(f(x),f(x'))\leq d_X(x,x')+\delta \qquad\text{for all }x,x'\in X,\\
				\sup_{y\in Y}d_Y(y,f(X))\leq \delta,\qquad \H_X^N(X)\leq\H_Y^N(Y)+\delta,\\
				\H_X^{N-1}(\partial X)\leq\H_Y^{N-1}(\partial Y)+\delta.
			\end{gathered}
		\end{equation}
		Then $f$ is an $\varepsilon$-Gromov--Hausdorff approximation, i.e.
		\begin{equation}
			Y\subseteq B_\varepsilon(f(X)),\quad\left|d_Y(f(x),f(x'))-d_X(x,x')\right|\leq\varepsilon
			\qquad\text{for all }x,x'\in X.
		\end{equation}
	\end{theorem}
	
	\begin{proof}
		Suppose by contradiction that there exist $\varepsilon_0>0$, ${\rm ncRCD}(K,N)$ spaces $X_i,Y_i$, and maps $f_i\colon X_i\to Y_i$ satisfying
		\eqref{eq:assumptions_almost_LVR} and \eqref{eq:defect_almost_LVR} with $\delta_i\downarrow0$, while $f_i$ is not an $\varepsilon_0$-Gromov--Hausdorff
		approximation.
		
		From the assumption \eqref{eq:defect_almost_LVR}, $\diam (Y_i)\leq 2D$ for all sufficiently large $i$.
		Thus all spaces considered are compact and have uniformly bounded diameters.
		Gromov's compactness theorem yields, up to passing to subsequences,
		\begin{equation}
			(X_i,d_{X_i})\xrightarrow{\mathrm{GH}}(X,d_X),
			\qquad
			(Y_i,d_{Y_i})\xrightarrow{\mathrm{GH}}(Y,d_Y),
		\end{equation}
		where $X,Y$ are compact metric spaces.
		
		Set
		$\Gamma_i\coloneqq\overline{\{(x,f_i(x)):x\in X_i\}}$. 
		Up to passing to a further subsequence, $\Gamma_i$ converges in the Hausdorff topology to a compact set $\Gamma\subseteq X\times Y$, whose projection onto $X$ is surjective.
		For any $(x,y),(x',y')\in\Gamma$, approximating points of the closed graphs by graph points, choose $x_i,x_i'\in X_i$ such that
		\begin{equation}
			x_i\to x,\qquad x_i'\to x',\qquad f_i(x_i)\to y,\qquad f_i(x_i')\to y'.
		\end{equation}
		Then
		\begin{equation}
			d_Y(y,y')=\lim_{i\to\infty}d_{Y_i}(f_i(x_i),f_i(x_i'))\leq\lim_{i\to\infty}\left(d_{X_i}(x_i,x_i')+\delta_i\right)= d_X(x,x').
		\end{equation}
		Taking $x=x'$ shows that $\Gamma$ is the graph of a map $f\colon X\to Y$, and the same inequality shows that $f$ is $1$-Lipschitz.
	    Moreover, graph convergence implies that $f_i(x_i)\to f(x)$ whenever $x_i\to x$.
		
		Surjectivity of $f$ follows from $Y_i\subseteq B_{2\delta_i}(f_i(X_i))$.
		Indeed, for any $y\in Y$, choose $y_i\in Y_i$ converging to $y$
		and $x_i\in X_i$ such that
		$d_{Y_i}(y_i,f_i(x_i))<2\delta_i$.
		Up to passing to a subsequence, $x_i\to x\in X$,
		and graph convergence gives $f(x)=y$.
		
	With the uniform volume lower bound \eqref{eq:assumptions_almost_LVR}, \cite[Theorem~1.2]{DePhilippis-Gigli18} ensures that
	$(Y,d_Y,\H_Y^N)$ is a ${\rm ncRCD}(K,N)$ space and the convergence of $(Y_i)$ is also measured Gromov--Hausdorff convergence.
	In particular, $\H_Y^N(Y)\geq v$.
	Since $f$ is surjective and $1$-Lipschitz, $\H^N(X)\geq\H_Y^N(Y)\geq v>0$.
	Again by \cite[Theorem~1.2]{DePhilippis-Gigli18}, $X$ cannot be a collapsed limit and consequently $(X,d_X,\H_X^N)$ is also a ${\rm ncRCD}(K,N)$ space.
	As a summary,
	\begin{equation}\label{eq:proof_mGH}
		(X_i,d_{X_i},\H_{X_i}^N)\xrightarrow{\mathrm{mGH}}(X,d_X,\H_X^N),\qquad
		(Y_i,d_{Y_i},\H_{Y_i}^N)\xrightarrow{\mathrm{mGH}}(Y,d_Y,\H_Y^N).
	\end{equation}
	From the assumption and the volume convergence,
	\begin{equation}
		\H_X^N(X)=\lim_{i\to\infty}\H_{X_i}^N(X_i)\leq\lim_{i\to\infty} \big(\H_{Y_i}^N(Y_i)+\delta_i\big)=\H_Y^N(Y).
	\end{equation} 
     Thus $\H_X^N(X)=\H_Y^N(Y)<\infty$. 
	Lemma~\ref{lemma:volume-preserving} ensures that $f$ is volume-preserving.
	Furthermore, applying the boundary measure convergence \cite[Theorem~1.8]{Brue-Naber-Semola22} to \eqref{eq:proof_mGH} with the assumption \eqref{eq:defect_almost_LVR}
	\begin{equation}
			\H_X^{N-1}(\partial X)=\lim_{i\to\infty}\H_{X_i}^{N-1}(\partial X_i)\leq\lim_{i\to\infty}
			\bigl(\H_{Y_i}^{N-1}(\partial Y_i)+\delta_i\bigr)=\H_Y^{N-1}(\partial Y).
	\end{equation}
	By \cref{prop:LVR_RCD_3}, $f$ is an isometry.
         
	For all sufficiently large $i$, $\delta_i<\varepsilon_0$.
		Hence the failure of the $\varepsilon_0$-Gromov--Hausdorff approximation property gives
		$x_i,x_i'\in X_i$
		\begin{equation}
			d_{X_i}(x_i,x_i')-d_{Y_i}(f_i(x_i),f_i(x_i'))\geq\varepsilon_0.
		\end{equation}
		Passing to a subsequence and then to the limit contradicts that $f$ is an isometry.
	\end{proof}

	\appendix
	\section{Remarks on the {\normalfont lip-to-Lip} property}\label{app:lip-to-Lip}
	In this appendix, we investigate more closely the lip-to-Lip property introduced in Definition \ref{def:L},
	discussing both necessary and sufficient conditions for its validity.
	\medskip
	
	We begin with some preliminary results concerning metric spaces.
	\begin{remark}\label{rmk:nec_cond_length}
		If \((X,d)\) is a length space, then it is well known that
		\begin{equation}\label{eq:Lip_sup_lip}
			{\rm Lip}(f)=\sup_{x\in X}{\rm lip}(f)(x)\quad\text{ for every Lipschitz function }f\colon X\to\R.
		\end{equation}
		Indeed, the inequality \(\geq\) follows from the trivial observation that (in an arbitrary metric space)
		\({\rm lip}(f)(x)\leq{\rm Lip}(f)\) holds for every \(x\in X\). To prove the converse inequality, fix any
		\(\varepsilon>0\) and pick points \(x,y\in X\), \(x\neq y\) so that \(\frac{f(y)-f(x)}{d(y,x)}\geq{\rm Lip}(f)-\varepsilon\).
		Take a constant-speed Lipschitz curve \(\gamma\colon[0,1]\to X\) such that \(\gamma_0=x\), \(\gamma_1=y\) and
		\(L_X(\gamma)<(1+\varepsilon)d(x,y)\). Hence,
		\[\begin{split}
			\int_0^1\lip(f)(\gamma_t)\,\d t&\geq\frac{1}{L_X(\gamma)}\int_0^1(f\circ\gamma)'(t)\,\d t
			=\frac{f(y)-f(x)}{L_X(\gamma)}>\frac{1}{1+\varepsilon}\frac{f(y)-f(x)}{d(y,x)}\\
			&\geq\frac{{\rm Lip}(f)}{1+\varepsilon}-\frac{\varepsilon}{1+\varepsilon},
		\end{split}\]
		whence it follows from \(\sup_{z\in X}\lip(f)(z)\ge \lip(f)(\gamma_{t_0})>(1+\varepsilon)^{-1}{\rm Lip}(f)-\varepsilon/(1+\varepsilon)\)
		for some \(t_0\in[0,1]\). By the arbitrariness of \(\varepsilon>0\), the inequality \(\leq\) in \eqref{eq:Lip_sup_lip}
		is established.
	\end{remark}

	The \textbf{asymptotic Lipschitz constant} function \(\lip_a(f)\colon X\to[0,+\infty)\) of a
	Lipschitz function \(f\colon X\to\R\) is defined as
	\[
	\lip_a(f)(x)\coloneqq\inf_{r>0}{\rm Lip}(f|_{B_r(x)})\quad\text{ for every }x\in X.
	\]
	It clearly holds that \(\lip(f)(x)\leq\lip_a(f)(x)\leq{\rm Lip}(f)\) for every \(x\in X\).
	\medskip
	
	In the next results, which we have not been able to find in the literature, we show that -- in the class of complete
	metric spaces -- property \eqref{eq:Lip_sup_lip} actually characterises length spaces.
	\begin{lemma}\label{lem:aux_equiv_compl_length_space}
		Let \((X,d)\) be a metric space. Assume that
		\begin{equation}\label{eq:aux_equiv_compl_length_space_hp}
			{\rm Lip}(f)=\sup_{x\in X}\lip_a(f)(x)\quad\text{ for every Lipschitz function }f\colon X\to\R.
		\end{equation}
		Then the metric completion \((\bar X,\bar d)\) of \((X,d)\) is a length space.
	\end{lemma}
	\begin{proof}
		Given a metric space \((Y,\rho)\) and $\varepsilon>0$, we define the $\varepsilon$-chain distance on \(Y\) as
		\[
		\rho_\varepsilon(x,y)\coloneqq
		\inf\left\{
		\sum_{i=1}^N\rho(y_{i-1},y_i)\;\Bigg|\;
		\begin{array}{l}
			y_0,\ldots,y_N\in Y,\,y_0=x,\,y_N=y,\\
			\rho(y_{i-1},y_i)<\varepsilon\text{ for every }i=1,\ldots,N
		\end{array}
		\right\}\quad\forall x,y\in Y.
		\]
		First, we claim that \(d_\varepsilon=d\) for all \(\varepsilon>0\). Clearly, \(d_\varepsilon\geq d\).
		To prove the converse inequality, fix $x\in X$ and $R>0$,
		and define $\phi_{x,R}(z)\coloneqq d_\varepsilon(x,z)\wedge R$ for all \(z\in X\). Note that
		\begin{equation}
			d_\varepsilon(x,w)\leq d_\varepsilon(x,z)+d(z,w)\quad\text{ for every }z,w\in X\text{ with }d(z,w)<\varepsilon.
		\end{equation}
		It follows that $\lip_a (\phi_{x,R})(z)\leq 1$ for all $z\in X$. Moreover, $\phi_{x,R}$ is globally
		Lipschitz with $\Lip(\phi_{x,R})$ bounded by $\max\{1,\frac{R}{\varepsilon}\}$, since for every
		\(z,w\in X\) with $d(z,w)\geq\varepsilon$ we can estimate
		\[
		|\phi_{x,R}(z)-\phi_{x,R}(w)|\leq R\leq \frac{R}{\varepsilon}d(z,w).
		\]
		Then \eqref{eq:aux_equiv_compl_length_space_hp} implies that $\phi_{x,R}$ is $1$-Lipschitz.
		As $\phi_{x,R}(x)=0$, we infer that $\phi_{x,R}(y)\leq d(x,y)$ for every \(y\in X\). Letting
		$R\to \infty$, we get $d_{\varepsilon}(x,y)\leq d(x,y)$ for every \(x,y\in X\).
		
		Next, we claim that \(\bar d_\varepsilon=\bar d\) for all \(\varepsilon>0\). By the $\varepsilon$-chain
		characterization of complete length spaces (see e.g.\ \cite[Corollary 2.4.17]{Burago-Burago-Ivanov}),
		this implies that $(\bar X,\bar d)$ is a length space. The inequality \(\bar d_\varepsilon\geq\bar d\)
		holds by the triangle inequality. To prove the converse inequality, fix \(x,y\in\bar X\). Given any \(\eta>0\),
		we take \(x',y'\in X\) such that \(\bar d(x,x'),\bar d(y,y')<\min\{\varepsilon,\frac{\eta}{8}\}\).
		We can find an \(\varepsilon\)-chain in \(X\), joining \(x'\) and \(y'\), whose total length is less
		that \(d(x',y')+\frac{\eta}{2}\). Extending this chain by adding the points \(x\) and \(y\), we obtain
		an $\varepsilon$-chain in $\bar X$ whose total length is less than $\bar d(x,y)+\eta$.
		Letting \(\eta\to 0\), we conclude that \(\bar d_\varepsilon(x,y)\leq\bar d(x,y)\).
	\end{proof}

	Note that \eqref{eq:aux_equiv_compl_length_space_hp} does not guarantee that \((X,d)\) itself is a length
	space: consider e.g.\ the punctured real line \(\R\setminus\{0\}\) equipped with the Euclidean distance.
	\begin{corollary}\label{cor:equiv_compl_length_space}
		Let \((X,d)\) be a complete metric space. Then the following conditions are equivalent:
		\begin{itemize}
			\item[\(\rm i)\)] \((X,d)\) is a length space.
			\item[\(\rm ii)\)] \({\rm Lip}(f)=\sup_{x\in X}\lip(f)(x)\) for every Lipschitz function \(f\colon X\to\R\).
			\item[\(\rm iii)\)] \({\rm Lip}(f)=\sup_{x\in X}\lip_a(f)(x)\) for every Lipschitz function \(f\colon X\to\R\).
		\end{itemize}
	\end{corollary}
	\begin{proof}
		The implication \({\rm i)}\Rightarrow{\rm ii)}\) was observed in Remark \ref{rmk:nec_cond_length}.
		The implication \({\rm ii)}\Rightarrow{\rm iii)}\) follows from the fact that \(\lip(f)\leq\lip_a(f)\leq\Lip(f)\).
		The implication \({\rm iii)}\Rightarrow{\rm i)}\) is a direct consequence of Lemma \ref{lem:aux_equiv_compl_length_space}.
	\end{proof}

	We now focus on metric measure spaces satisfying the lip-to-Lip property.
	\begin{proposition}\label{prop:(L)-1}
		Let $(X,d,\m)$ be a metric measure space having the lip-to-Lip property.
		Then $\supp (\m)=X$. Moreover, the metric completion $(\bar X,\bar d)$ of $(X,d)$ is a length space.
	\end{proposition}
	\begin{proof}
		To prove that \(\m\) is fully supported, we argue by contradiction: suppose there exists \(x_0\in X\setminus\supp(\m)\).
		Then we can find \(r>0\) such that \(\m(B_{2r}(x_0))=0\). Consider the Lipschitz function \(\phi\), which we define
		as $\phi(x)\coloneqq(r-d(x,x_0))_+$ for every \(x\in X\). Then we have \(\Lip(\phi)=1\), whereas \(\lip(\phi)=0\) holds
		\(\m\)-a.e.\ on \(X\), thus leading to a contradiction.
		
		The lip-to-Lip property implies that \eqref{eq:aux_equiv_compl_length_space_hp} holds, since
		\(\|\lip(\phi)\|_{L^\infty(\m)}\leq\sup_{x\in X}\lip_a(\phi)(x)\). Therefore, Lemma \ref{lem:aux_equiv_compl_length_space}
		ensures that \((\bar X,\bar d)\) is a length space.
	\end{proof}

	On the contrary, a fully-supported measure on a length space does not necessarily yield a metric measure space satisfying
	the lip-to-Lip property, as the following example illustrates.
	\begin{example}
		Enumerate the rational numbers \(\mathbb Q\) as \(\{q_n:n\in\N\}\), and consider the fully-supported probability measure
		\(\m\coloneqq\sum_{n=1}^\infty 2^{-n}\delta_{q_n}\) on \(\R\). Then we claim that \((\R,d_{\rm Eucl},\m)\) does not have
		the lip-to-Lip property. To prove it, denote by \(g\) the characteristic function of \(C\coloneqq\R\setminus\bigcup_{n\in\N}(q_n-2^{-n},q_n+2^{-n})\),
		and we define \(f\in\Lip(\R)\) as
		\[
		f(t)\coloneqq\int_0^t g(s)\,\d s\quad\text{ for every }t\in\R.
		\]
		Since \(\mathcal L^1(\R\setminus C)\leq\sum_{n=1}^\infty 2^{-n+1}<+\infty\), we have that
		\[
		\frac{f(t)-f(-t)}{t-(-t)}=\frac{1}{2t}\int_{-t}^t g(s)\,\d s=\frac{\mathcal L^1((-t,t)\cap C)}{2t}\to 1\quad\text{ as }t\to+\infty,
		\]
		so that \(\Lip(f)=1\). Moreover, \(\lip(f)(q_n)=0\) for all \(n\in\N\), thus \(\|\lip(f)\|_{L^\infty(\m)}=0\).
	\end{example}
	Let us now consider sufficient conditions for the validity of the lip-to-Lip property.
	\begin{definition}
		A metric measure space $(X,d,\m)$ is called an \textbf{essential length space} if for all $x,y\in X$, all $N\subset X$ with $\m(N)=0$ and all $\varepsilon>0$, there is a Lipschitz curve $\gamma\colon [0,1]\to X$ connecting $x$ and $y$ such that $\L^1(\gamma^{-1}(N))=0$ and 
		\[
		d(x,y)+\varepsilon\geq L_X(\gamma).
		\]
	\end{definition}
	
	\begin{lemma}
		Every essentially length space $(X,d,\m)$ satisfies the lip-to-Lip property.
	\end{lemma}
	\begin{proof}
		Fix \(\phi\in\Lip(X)\) and set \(L\coloneqq\|\lip(\phi)\|_{L^\infty(\m)}\). We aim to show that \(\phi\) is \(L\)-Lipschitz,
		whence it follows that \(\Lip(\phi)=L\). Note that the set $N\coloneqq\{x\in X:\lip(\phi)(x)>L\}$ is $\m$-negligible. Hence, given
		any \(x,y\in X\) and \(\varepsilon>0\), we can find a Lipschitz curve \(\gamma\colon[0,1]\to X\) joining \(x\) and \(y\) such that
		$\L^1(\gamma^{-1}(N))=0$ and \(d(x,y)+\varepsilon\geq L_X(\gamma)\). Denoting by $|\dot{\gamma}_t|$ the metric derivative of $\gamma$
		and using the Lipschitzianity of \(\phi\circ\gamma\), we can thus estimate
		\[
		\lip(\phi\circ\gamma)(t)\leq\lip(\phi)(\gamma_t)|\dot\gamma_t|\leq L|\dot\gamma_t|\quad\text{ for }\L^1\text{-a.e.\ }t\in[0,1].
		\]
		It follows that
		\[
		|\phi(y)-\phi(x)|=|(\phi\circ\gamma)(1)-(\phi\circ\gamma)(0)|\leq\int_0^1\lip(\phi\circ\gamma)(t)\,\d t\leq L\int_0^1|\dot\gamma_t|\,\d t
		\leq L(d(x,y)+\varepsilon).
		\]
		Letting \(\varepsilon\to 0\), we conclude that \(|\phi(y)-\phi(x)|\leq L\,d(x,y)\) for every \(x,y\in X\), as desired.
	\end{proof}
	\begin{lemma}\label{lemma:essentiallength}
		Let $(X,d,\m)$ be a metric measure space satisfying the following condition: for any $x\in X$ and any Borel $A\subset X$ with $\m(A)\in (0,\infty)$, there exists a Borel probability measure $\pi$ on the space of geodesics ${\rm Geo}([0,1],X)\subset C([0,1],X)$ such that
		\[
		(e_0)_\#\pi=\delta_x,\quad(e_1)_{\#}\pi=\frac{\m|_A}{\m(A)},\quad(e_t)_{\#}\pi\ll \m \;\;\forall t\in(0,1],
		\]
		where $e_t\colon \gamma\mapsto \gamma_t$ is the evaluation map.
		Then $(\supp(\m),d,\m)$ is an essential length space.
	\end{lemma}
	\begin{proof}
		Fix \(x,y\in\supp(\m)\), an \(\m\)-negligible set \(N\subset X\), and \(\varepsilon>0\).
		Applying the assumption with \(A=B_\varepsilon(y)\), we obtain Borel probabilities $\pi^x,\pi^y$ on $\rm{Geo}([0,1],X)$ such that
		\[
		(e_0)_\#\pi^x=\delta_x,\quad (e_1)_\#\pi^y=\delta_y,\quad (e_1)_\#\pi^x=(e_0)_\#\pi^y=\frac{\m|_A}{\m(A)}
		\]
		and $(e_t)_\#\pi^p\ll\m$ for every $t\in (0,1)$ and $p\in\{x,y\}$. Let us define \(G\coloneqq B_\varepsilon(y)\setminus(B_x\cup B_y)\), where we set
		\[\begin{split}
			B_x&\coloneqq\big\{z\in B_\varepsilon(y)\;\big|\;\L^1(\gamma^{-1}(N))>0\text{ for every }\gamma\in\Geo([0,1],X)\cap e_0^{-1}(x)\cap e_1^{-1}(z)\big\},\\
			B_y&\coloneqq\big\{z\in B_\varepsilon(y)\;\big|\;\L^1(\sigma^{-1}(N))>0\text{ for every }\sigma\in\Geo([0,1],X)\cap e_0^{-1}(z)\cap e_1^{-1}(y)\big\},
		\end{split}\]
		By Fubini's theorem, we have that
		\begin{align}
			\int_{\Geo([0,1],X)}\int^1_0 \mathds{1}_{N}(\gamma_t)\d t\d \pi^p(\gamma)=\int^1_0 (e_t)_{\#}\pi^p(N)\d t=0\quad\text{ for every }p\in\{x,y\},
		\end{align}
		so that \(\L^1(\gamma^{-1}(N))=0\) for \(\pi^x\)-a.e.\ \(\gamma\) and \(\L^1(\sigma^{-1}(N))=0\) for \(\pi^y\)-a.e.\ \(\sigma\).
		In particular, we have that \(\pi^x(e_1^{-1}(B_x))=0\) and thus \(\m(B_x)=\m(B_\varepsilon(y))(e_1)_\#\pi^x(B_x)=0\), and similarly
		\(\m(B_y)=0\). This gives that \(\m(G)>0\), so that \(G\neq\varnothing\). Now, fix any \(z\in G\). Then we can find geodesics
		\(\gamma,\sigma\colon[0,1]\to X\) such that \(\gamma_0=x\), \(\gamma_1=z=\sigma_0\), \(\sigma_1=y\) and \(\L^1(\gamma^{-1}(N))=\L^1(\sigma^{-1}(N))=0\). Therefore, the curve \(\gamma*\sigma\) obtained by concatenating \(\gamma\) and \(\sigma\) is a Lipschitz curve joining
		\(x\) and \(y\) that satisfies \(\L^1((\gamma*\sigma)^{-1}(N))=0\) and
		\[
		L_X(\gamma*\sigma)=L_X(\gamma)+L_X(\sigma)=d(x,z)+d(z,y)\leq d(x,y)+2\varepsilon.
		\]
		This shows that $(\supp(\m),d,\m)$ is an essential length space, as desired.
	\end{proof}
	
	In particular, \cref{lemma:essentiallength} applies when $(X,d,\m)$ satisfies the ${\rm MCP}$ condition in the sense of \cite{Ohta-MCP07,Sturm06-2}, or the strong interpolation property in \cite{Kell17}.
	We can summarize the above discussions as follows (with the last implication under the completeness assumption):
	\[
	\text{MCP space }\Rightarrow\text{ essential length space }\Rightarrow\text{ lip-to-Lip property }\Rightarrow\text{ length space.}
	\]
	Furthermore, another sufficient condition for the validity of the lip-to-Lip property is the so-called
	\emph{\(p\)-Sobolev-to-Lipschitz property}, considered in \cite{AmbrosioGigliSavare2015BakryEmery} (see also
	\cite{DelloSchiavoSuzuki2022SobolevToLipschitz,CreutzSoultanis2020MaximalMetricSurfaces,CaputoCavallucciIkonen2026SobolevToLipschitz} for further developments).
	We recall that a metric measure space \((X,d,\m)\) has the \(p\)-Sobolev-to-Lipschitz property, for some exponent \(p\in[1,\infty]\),
	if the following holds: every Newtonian function \(f\in N^{1,p}(X)\) whose minimal \(p\)-weak upper gradient \(g_f\)
	belongs to \(L^\infty(\m)\) has a Lipschitz representative \(\bar f\) satisfying \(\Lip(\bar f)=\|g_f\|_{L^\infty(\m)}\).
	We refer e.g.\ to the monograph \cite{HeinonenKoskelaShanmugalingamTyson2015} for more details on the Newtonian space.
	\begin{lemma}\label{lem:Sob-to-Lip_implies_(L)}
		Let \((X,d,\mathfrak m)\) be a metric measure space with \(\supp(\m)=X\) having the \(p\)-Sobolev-to-Lipschitz property for some \(p\in[1,\infty]\).
		Then \((X,d,\mathfrak m)\) has the lip-to-Lip property.
	\end{lemma}
	\begin{proof}
		By a standard cut-off argument, it suffices to check that \({\rm Lip}(f)=\|{\rm lip}(f)\|_{L^\infty(\mathfrak m)}\)
		for every boundedly-supported Lipschitz function \(f\colon X\to\R\). Since we have that \(f\in N^{1,p}(X)\) and
		\(g_f\leq{\rm lip}(f)\leq\|{\rm lip}(f)\|_{L^\infty(\mathfrak m)}\)
		holds \(\mathfrak m\)-a.e.\ on \(X\), the \(p\)-Sobolev-to-Lipschitz property
		yields \({\rm Lip}(f)=\|g_f\|_{L^\infty(\m)}\leq\|{\rm lip}(f)\|_{L^\infty(\mathfrak m)}\).
		The converse inequality trivially holds.
	\end{proof}

	Combining Lemma \ref{lem:Sob-to-Lip_implies_(L)} with the results of \cite{AmbrosioGigliSavare2015BakryEmery},
	we deduce that every \({\rm RCD}(K,\infty)\) space has the lip-to-Lip property.
	%
	%
	%
	%

\begin{thebibliography}{CDNZ{\etalchar{+}}24}

\bibitem[AGS14]{AGS14-RCD}
Luigi Ambrosio, Nicola Gigli, and Giuseppe Savar\'{e}.
\newblock Metric measure spaces with {R}iemannian {R}icci curvature bounded from below.
\newblock {\em Duke Math. J.}, 163(7):1405--1490, 2014.

\bibitem[AGS15]{AmbrosioGigliSavare2015BakryEmery}
Luigi Ambrosio, Nicola Gigli, and Giuseppe Savar{\'e}.
\newblock Bakry--\'{E}mery curvature-dimension condition and riemannian {Ricci} curvature bounds.
\newblock {\em The Annals of Probability}, 43(1):339--404, 2015.

\bibitem[AK00]{Ambrosio-Kirchheim00}
Luigi Ambrosio and Bernd Kirchheim.
\newblock Rectifiable sets in metric and {Banach} spaces.
\newblock {\em Math. Ann.}, 318(3):527--555, 2000.

\bibitem[BBI01]{Burago-Burago-Ivanov}
D.~Burago, Yu. Burago, and S.~Ivanov.
\newblock {\em A course in metric geometry}, volume~33 of {\em Grad. Stud. Math.}
\newblock Providence, RI: American Mathematical Society (AMS), 2001.

\bibitem[BCS23]{BassoCreutzSoultanis2023}
Giuliano Basso, Paul Creutz, and Elefterios Soultanis.
\newblock Filling minimality and {L}ipschitz-volume rigidity of convex bodies among integral current spaces.
\newblock {\em Journal f{\"u}r die reine und angewandte Mathematik}, 805:213--239, 2023.

\bibitem[BH99]{Bridson-Haefliger99}
Martin~R. Bridson and Andr{\'e} Haefliger.
\newblock {\em Metric Spaces of Non-Positive Curvature}, volume 319 of {\em Grundlehren der mathematischen Wissenschaften}.
\newblock Springer-Verlag, Berlin, 1999.

\bibitem[BNS22]{Brue-Naber-Semola22}
Elia Bru{\'e}, Aaron Naber, and Daniele Semola.
\newblock Boundary regularity and stability for spaces with {Ricci} bounded below.
\newblock {\em Invent. Math.}, 228(2):777--891, 2022.

\bibitem[BT24]{Bate-Takavc24}
David Bate and Jakub Tak{\'a}{\v{c}}.
\newblock Typical {Lipschitz} images of rectifiable metric spaces.
\newblock {\em J. Reine Angew. Math.}, 810:139--188, 2024.

\bibitem[CC97]{Cheeger-Colding97}
Jeff Cheeger and Tobias~H. Colding.
\newblock On the structure of spaces with {Ricci} curvature bounded below. {I}.
\newblock {\em J. Differ. Geom.}, 46(3):406--480, 1997.

\bibitem[CCI26]{CaputoCavallucciIkonen2026SobolevToLipschitz}
Emanuele Caputo, Nicola Cavallucci, and Toni Ikonen.
\newblock Sobolev-to-{Lipschitz} property of geodesically complete spaces with curvature bounded from above.
\newblock Preprint, arXiv:2606.24374, 2026.

\bibitem[CDNZ{\etalchar{+}}24]{CDNPSG24}
Chris Connell, Xianzhe Dai, Jes{\'u}s N{\'u}{\~n}ez-Zimbr{\'o}n, Raquel Perales, Pablo Su{\'a}rez-Serrato, and Guofang Wei.
\newblock Volume entropy and rigidity for {RCD}-spaces.
\newblock Preprint, {arXiv}:2411.04327 [math.{DG}] (2024), 2024.

\bibitem[CS20]{CreutzSoultanis2020MaximalMetricSurfaces}
Paul Creutz and Elefterios Soultanis.
\newblock Maximal metric surfaces and the {S}obolev-to-{Lipschitz} property.
\newblock {\em Calculus of Variations and Partial Differential Equations}, 59(5):177, 2020.

\bibitem[Den25]{Deng25}
Qin Deng.
\newblock H{\"o}lder continuity of tangent cones in {RCD}{$(K,N)$} spaces and applications to nonbranching.
\newblock {\em Geom. Topol.}, 29(2):1037--1114, 2025.

\bibitem[DNP23]{DelNinPerales2023}
Giacomo Del~Nin and Raquel Perales.
\newblock Rigidity of mass-preserving {1}-{L}ipschitz maps from integral current spaces into {$\mathbb{R}^n$}.
\newblock {\em Journal of Mathematical Analysis and Applications}, 526(1):127297, 2023.

\bibitem[DPG18]{DePhilippis-Gigli18}
Guido De~Philippis and Nicola Gigli.
\newblock Non-collapsed spaces with {Ricci} curvature bounded from below.
\newblock {\em J. {\'E}c. Polytech., Math.}, 5:613--650, 2018.

\bibitem[DPMR17]{DePhilippisMarcheseRindler2017}
Guido De~Philippis, Andrea Marchese, and Filip Rindler.
\newblock On a conjecture of {Cheeger}.
\newblock In Nicola Gigli, editor, {\em Measure Theory in Non-Smooth Spaces}, pages 145--155. De Gruyter Open, 2017.

\bibitem[DSS22]{DelloSchiavoSuzuki2022SobolevToLipschitz}
Lorenzo Dello~Schiavo and Kohei Suzuki.
\newblock Sobolev-to-{Lipschitz} property on {QCD}-spaces and applications.
\newblock {\em Mathematische Annalen}, 384(3--4):1815--1832, 2022.

\bibitem[Fed69]{Federer69}
Herbert Federer.
\newblock {\em Geometric measure theory}, volume 153 of {\em Grundlehren Math. Wiss.}
\newblock Springer, Cham, 1969.

\bibitem[Gig15]{Gigli15-RCD}
Nicola Gigli.
\newblock {\em On the differential structure of metric measure spaces and applications}, volume 1113 of {\em Mem. Am. Math. Soc.}
\newblock Providence, RI: American Mathematical Society (AMS), 2015.

\bibitem[GP20]{Gigli-Pasqualetto-Book}
Nicola Gigli and Enrico Pasqualetto.
\newblock {\em Lectures on nonsmooth differential geometry}, volume~2 of {\em SISSA Springer Ser.}
\newblock Cham: Springer, 2020.

\bibitem[GP21]{GigliPasqualetto2021}
Nicola Gigli and Enrico Pasqualetto.
\newblock Behaviour of the reference measure on {RCD} spaces under charts.
\newblock {\em Communications in Analysis and Geometry}, 29(6):1391--1414, 2021.

\bibitem[GP22]{Gigli-Pasqualetto22}
Nicola Gigli and Enrico Pasqualetto.
\newblock Equivalence of two different notions of tangent bundle on rectifiable metric measure spaces.
\newblock {\em Commun. Anal. Geom.}, 30(1):1--51, 2022.

\bibitem[GRS16]{Gigli-Rajala-Sturm16}
Nicola Gigli, Tapio Rajala, and Karl-Theodor Sturm.
\newblock Optimal maps and exponentiation on finite-dimensional spaces with {Ricci} curvature bounded from below.
\newblock {\em J. Geom. Anal.}, 26(4):2914--2929, 2016.

\bibitem[HKST15]{HeinonenKoskelaShanmugalingamTyson2015}
Juha Heinonen, Pekka Koskela, Nageswari Shanmugalingam, and Jeremy~T. Tyson.
\newblock {\em Sobolev Spaces on Metric Measure Spaces: An Approach Based on Upper Gradients}, volume~27 of {\em New Mathematical Monographs}.
\newblock Cambridge University Press, 2015.

\bibitem[Kel17]{Kell17}
Martin Kell.
\newblock Transport maps, non-branching sets of geodesics and measure rigidity.
\newblock {\em Adv. Math.}, 320:520--573, 2017.

\bibitem[Kir94]{Kirchheim94}
Bernd Kirchheim.
\newblock Rectifiable metric spaces: {Local} structure and regularity of the {Hausdorff} measure.
\newblock {\em Proc. Am. Math. Soc.}, 121(1):113--123, 1994.

\bibitem[KL16]{Kitabeppu-Lakzian16}
Yu~Kitabeppu and Sajjad Lakzian.
\newblock Characterization of low dimensional {{\(RCD^\ast(K, N)\)}} spaces.
\newblock {\em Anal. Geom. Metr. Spaces}, 4:187--215, 2016.

\bibitem[KM18]{KellMondino2018}
Martin Kell and Andrea Mondino.
\newblock On the volume measure of non-smooth spaces with {Ricci} curvature bounded below.
\newblock {\em Annali della Scuola Normale Superiore di Pisa, Classe di Scienze}, 18(2):593--610, 2018.

\bibitem[KM21]{Kapovitch-Mondino21}
Vitali Kapovitch and Andrea Mondino.
\newblock On the topology and the boundary of {{\(N\)}}-dimensional {{\(\mathsf{RCD}\,(K,N)\)}} spaces.
\newblock {\em Geom. Topol.}, 25(1):445--495, 2021.

\bibitem[Li15]{Li15}
Nan Li.
\newblock Lipschitz-volume rigidity in {Alexandrov} geometry.
\newblock {\em Adv. Math.}, 275:114--146, 2015.

\bibitem[LS23]{Lytchak-Stadler23}
Alexander Lytchak and Stephan Stadler.
\newblock Ricci curvature in dimension 2.
\newblock {\em J. Eur. Math. Soc. (JEMS)}, 25(3):845--867, 2023.

\bibitem[LV09]{LottVillani09}
John Lott and Cedric Villani.
\newblock Ricci curvature for metric-measure spaces via optimal transport.
\newblock {\em Ann. Math. (2)}, 169(3):903--991, 2009.

\bibitem[LW14]{Li-Wang14}
Nan Li and Feng Wang.
\newblock Lipschitz-volume rigidity on limit spaces with {Ricci} curvature bounded from below.
\newblock {\em Differ. Geom. Appl.}, 35:50--55, 2014.

\bibitem[Mar26]{Marti26}
Denis Marti.
\newblock The {L}ipschitz-volume rigidity problem for metric manifolds.
\newblock {\em The Journal of Geometric Analysis}, 36, 10 2026.

\bibitem[MMPZ03]{Mera-Moran-Preiss03}
M.~E. Mera, M.~Mor{\'a}n, D.~Preiss, and L.~Zaj{\'{\i}}{\v{c}}ek.
\newblock Porosity, {{\(\sigma\)}}-porosity and measures.
\newblock {\em Nonlinearity}, 16(1):247--255, 2003.

\bibitem[MN19]{Mondino-NaberJEMS}
Andrea Mondino and Aaron Naber.
\newblock Structure theory of metric measure spaces with lower {Ricci} curvature bounds.
\newblock {\em J. Eur. Math. Soc. (JEMS)}, 21(6):1809--1854, 2019.

\bibitem[MN24]{Meier-Ntalaepekos24}
Damaris Meier and Dimitrios Ntalampekos.
\newblock Lipschitz-volume rigidity and {Sobolev} coarea inequality for metric surfaces.
\newblock {\em J. Geom. Anal.}, 34(5):30, 2024.
\newblock Id/No 128.

\bibitem[Oht07]{Ohta-MCP07}
Shin-ichi Ohta.
\newblock On the measure contraction property of metric measure spaces.
\newblock {\em Comment. Math. Helv.}, 82(4):805--828, 2007.

\bibitem[Stu06]{Sturm06-2}
Karl-Theodor Sturm.
\newblock On the geometry of metric measure spaces. {II}.
\newblock {\em Acta Math.}, 196(1):133--177, 2006.

\bibitem[Z{\"u}s24]{Zust24}
Roger Z{\"u}st.
\newblock Lipschitz-volume rigidity of {Lipschitz} manifolds among integral currents.
\newblock {\em J. Geom. Anal.}, 34(7):20, 2024.
\newblock Id/No 210.

\end{thebibliography}
%
%
\newcommand{\etalchar}[1]{$^{#1}$}

\end{document}